\documentclass{article}

\usepackage{PRIMEarxiv}
 \usepackage{amsmath}
\usepackage{amssymb, nicefrac}
\usepackage{color}

\usepackage[ruled,vlined]{algorithm2e}
\usepackage{diagbox}
\usepackage{subcaption}

\usepackage[utf8]{inputenc} % allow utf-8 input
\usepackage[T1]{fontenc}    % use 8-bit T1 fonts
\usepackage{hyperref}       % hyperlinks
\usepackage{url}            % simple URL typesetting
\usepackage{booktabs}       % professional-quality tables
\usepackage{amsfonts}       % blackboard math symbols
\usepackage{nicefrac}       % compact symbols for 1/2, etc.
\usepackage{microtype}      % microtypography
\usepackage{lipsum}
\usepackage{fancyhdr}       % header
\usepackage{graphicx}       % graphics
\graphicspath{{media/}}     % organize your images and other figures under media/ folder

\usepackage{endnotes}
\let\footnote=\endnote

\newtheorem{assumption}{Assumption}
\newtheorem{proposition}{Proposition}
\newtheorem{lemma}{Lemma}
\usepackage{natbib}
 \bibpunct[, ]{(}{)}{,}{a}{}{,}%

\newcommand{\bxi}{\mbox{\boldmath$\xi$}}

\newcommand{\var}{{\rm Var}}

\newcommand{\U}{{\cal U}}
\newcommand{\V}{{\cal V}}

\newcommand{\B}{{\cal B}}

\newcommand{\I}{{\cal I}}

\newcommand{\X}{{\cal X}}

\def\bbr{{\Bbb{R}}} %real numbers
\def\bbe{{\Bbb{E}}} %expectation
\def\bbn{{\Bbb{N}}}

\def\bbb{{\Bbb{B}}}

\def\dist{\mathop{\rm dist}}

\def\e{\varepsilon}

\def\vv{\vartheta}

\newcommand{\cT}{{\mathfrak T}}

\newcommand{\N}{{\cal N}}

\newcommand{\T}{{\cal T}}

\newcommand{\half}{ \mbox{\small$\frac{1}{2}$}}

\newcommand{\ind}{{\mbox{\boldmath$1$}}}

\title{Episodic Bayesian Optimal Control with Unknown Randomness Distributions 
}

\author{
  Alexander Shapiro, Enlu Zhou  \\
  School of Industrial and Systems Engineering \\
  Georgia Institute of Technology \\
  Atlanta\\
  \texttt{\{ashapiro, enlu.zhou\}@isye.gatech.edu} \\
   \And
  Yifan Lin, Yuhao Wang \\
  School of Industrial and Systems Engineering \\
  Georgia Institute of Technology \\
  Atlanta\\
  \texttt{\{yifanlin,yuhaowang\}@gatech.edu} \\
}

\begin{document}
\maketitle

\begin{abstract}
Stochastic optimal control with unknown randomness distributions has been studied for a long time, encompassing robust control, distributionally robust control, and adaptive control. We propose a new episodic Bayesian approach that incorporates Bayesian learning with optimal control. In each episode, the approach learns the randomness distribution with a Bayesian posterior and subsequently solves the corresponding Bayesian average estimate of the true problem. The resulting policy is exercised during the episode, while additional data/observations of the randomness are collected to update the Bayesian posterior for the next episode. We show that the resulting episodic value functions and policies converge almost surely to their optimal counterparts of the true problem  if the parametrized model of the randomness distribution is correctly specified. We further show that the asymptotic convergence rate of the episodic value functions is of the order $O(N^{-1/2})$, where $N$ is the number of episodes given that only one data
point is collected in each episode. We develop an efficient computational method based on stochastic dual dynamic programming (SDDP) for a class of problems that   have convex cost functions and linear state dynamics. Our numerical results on a classical inventory control problem verify the theoretical convergence results, and numerical comparison with two other methods demonstrate the effectiveness of the proposed Bayesian approach.
\end{abstract}

% keywords can be removed

\section{Introduction} \label{sec: intro}

Stochastic Optimal Control (SOC) provides a principled approach to dynamic decision making under uncertainty. It models the transition of the system state with a dynamic equation driven by randomness, with the assumption that the distribution of randomness is known. However, in practical problems, modeling randomness frequently depends on pre-collected data or observations, which introduces uncertainty regarding the distribution of randomness. This type of model uncertainty is also referred to as Knightian uncertainty or epistemic uncertainty in the literature.

There has been a long history of efforts to address this model uncertainty in stochastic optimal control, ranging from robust control to the recently developed distributionally robust control. Robust control (e.g., \cite{Gilboa1989MaxminEU}, \cite{Srbu2014ANO}, \cite{Hansen2006robust}), closely related to minimax control and $H^\infty$ control (e.g., \cite{Gonz2002MinimaxControl}, \cite{Basar2008Hinfty}), often models the control problem as a game between the decision maker and  nature, where nature selects the worst-case scenario in the uncertainty set while the decision maker aims to find the optimal control in the worst-case scenario. Depending on whether the uncertainty set remains fixed throughout the process or is constructed at each stage, the game between the decision maker and nature can be static or dynamic. In  recent years, the vibrant development of Distributionally Robust Optimization (DRO) has sparked research in applying this approach to stochastic optimal control. DRO is conceptually similar to robust optimization, but instead of an uncertainty set for unknown parameters, it constructs an ambiguity set for unknown distributions. For instance,  \cite{Parys2016DRC} applied DRO with an ambiguity set based on moment constraints to control systems with linear dynamics and quadratic costs; \cite{Tzortzis2019TVA} solved infinite horizon average cost control problems using total variation distance ambiguity sets on the conditional distribution of the controlled process; and \cite{Yang2018WassersteinDR} adopted a DRO approach with Wasserstein-based ambiguity sets for stochastic control problems. Both robust control and distributionally robust control have faced the criticism of being overly conservative, as they focus on worst-case scenarios that often rarely occur in practice. Furthermore, these approaches are typically non-adaptive, meaning their uncertainty or ambiguity sets are solely based on pre-collected data before the process begins, without incorporating observations of randomness revealed during the process. It is worth mentioning that there are also exceptions where works such as \cite{Lim2006Robust} and \cite{Bielecki2019ARC} incorporate learning into robust approaches, enabling the reduction of model uncertainty through observed data. In particular, a Bayesian risk optimization (BRO) formulation has been proposed in \cite{Wu:2018BRO} as an alternative to the DRO formulation, and it is more amenable to incorporate Bayesian learning for dynamic problems. For example, BRO has been extended to Markov decision processes by \cite{Zhou2022NIPS} and to data-driven stochastic optimization with streaming input data by \cite{ZhouLiu2022bayesianSGD}.

Regarding continuing learning from observations of randomness, the general approach of adaptive control incorporates model learning into optimization of the control policy and often takes into account the process of future revealed randomness. Bayesian learning, a natural way to learn from sequentially arrived data, is frequently used in adaptive control. The Bayesian posterior distribution is often treated as a state, sometimes referred to as the belief state or information state, whose dynamics are governed by the Bayes updating rule and driven by randomness. Dynamic programming is then employed on the augmented pair of the posterior (belief state) and the original physical state to obtain a policy (e.g., \cite{Kumar2015BayesianAC}, \cite{Rieder1975BayesianDP}). This approach essentially views the unknown distribution as a partially observed state and leads to the optimal policy in principle. However, it is often computationally prohibitive to solve, except for some special cases, due to the potentially infinite dimension of the Bayesian posterior. Even when a conjugate prior family of distribution is used, resulting in a closed-form and finite-dimensional posterior, the posterior remains a continuous state and hence solving the associate dynamic programming runs into the notorious curse of dimensionality. On a related note, Bayesian adaptive approaches have also been actively explored for Markov decision processes with parameter uncertainty (e.g., \cite{duff2002optimal}, \cite{Zhou2022NIPS}), and numerical methods have been developed to combat the aforementioned computational difficulty brought by the posterior (belief state). However, these methods  that treat posterior as a belief state are developed for finite-horizon problems, and their approximation error accrues as the horizon increases, making them inapplicable for the infinite-horizon problems such as the one considered in this paper.

In this paper, we consider an infinite-horizon discounted-total-cost SOC problem with an unknown randomness distribution. Similar to adaptive control, we adopt a parametric Bayesian approach to learn from observed data of the randomness over time. However, in view of computational challenges associated with treating the Bayesian posterior as a state, we take a suboptimal but computationally easier approach by considering the Bayesian average counterpart of the original problem. This formulation ignores the future revelation of the randomness process, and hence the resultant control policy is suboptimal compared to the approach of taking into account the dynamics of the posterior (belief state). Nevertheless, the resultant suboptimal policy is only exercised for one episode while more data are collected during this episode; these data are then used to update the Bayesian posterior and the Bayesian average problem for the next episode, and a new improved policy is computed for use in that episode. By continually updating Bayesian estimation with more observed data and improving the policy in such an iterative way, the estimated (Bayesian average) problem converges to the true problem and the resultant policy converges to the true optimal policy.

%{\color{blue} 
The most relevant approach in the literature is probably the posterior sampling approach for reinforcement learning (PSRL) which was first proposed in \cite{strens2000bayesian}. The regret analysis of PSRL with a fixed episode length (i.e. fixed frequency of updating the posterior distribution) has been studied in \cite{osband2013more,osband2014near,osband2017posterior,theocharous2017posterior,abeille2018improved}. A variant of PSRL, namely LazyPSRL that uses adaptive episodic length was studied in 
\cite{osband2013more,abbasi2015bayesian}. The main difference between the posterior sampling approach and our proposed Bayesian optimal control is how the posterior distribution is utilized. To be specific, PSRL draws one sample from the posterior and then solves the problem under this sampled environment, whereas our algorithm takes into account the entire posterior and solves the corresponding Bayesian average counterpart of the objective.   PSRL is developed for reinforcement learning, and hence the posterior sampling approach is useful for balancing the trade-off between exploration and exploitation. In contrast, we consider a control problem, and therefore, exploiting the information in the entire posterior is more beneficial. This distinction was also supported by our numerical results, which show that the algorithm performance improves as the sample size from the posterior increases. Moreover, we develop a cutting-plane-type computational method for continuous-state and continuous-action problems, alleviating the curse of dimensionality usually encountered by approaches that rely on some form of discretizing the state and action spaces. Theoretically, PSRL studies regret bounds of the episodic value function by making assumptions on the concentration of the Bayesian posterior, whereas we invoke the Bernstein-von Mises theorem (also known as the Bayesian central limit theorem) to derive the asymptotic convergence rate of the episodic value function. It will be an interesting future direction to extend our analysis approach to the PSRL stream of work.

In summary, our paper makes the following contributions.
\begin{itemize}
    \item We propose a novel formulation for stochastic optimal control problems with unknown randomness distributions. This formulation adapts to the observed randomness by incorporating Bayesian learning into control, while maintains computational feasibility  by fixing the Bayesian posterior in each episode. 

    \item  % Under regularity conditions, the Bayesian posterior converges in probability to the true parameter of the randomness distribution if our parametric model is correctly specified. Consequently, 
    We rigorously prove that the sequence of episodic value functions (i.e., the optimal value function of the Bayesian average problem in each episode) and associated policies  converge almost surely to the true optimal ones of the original problem, given that the distribution family of the randomness is correctly specified. Moreover, with an approximation commonly used in statistical analysis we derive the asymptotics of the episodic value function, which suggest its asymptotic convergence rate is of the order $O(N^{-1/2})$, where $N$ is the number of episodes if only one data point is collected in each episode. As a demonstration, we also analytically derive the asymptotic convergence rate on a classical inventory control problem and further verify it numerically. % the theoretical results in a numerical example of inventory control problem. 

    \item To solve our proposed formulation, we develop an algorithm based on the Stochastic Dual Dynamic Programming  (SDDP) method for a class of problems that have convex cost functions and linear state dynamics. The algorithm iteratively adds cutting planes to approximate the convex value function. To further improve the computation efficiency, we propose a warm-start of SDDP iterations in each episode by reusing the cutting planes from the previous episode.  We note that it is the first time that SDDP is extended to the infinite-horizon episodic setting. Due to the scarcity of existing convergence and complexity results of SDDP for the infinite-horizon models (theoretical properties of SDDP for finite-stage models can be found in a recent survey paper \cite{LanShapiro:2024survey}), rigorously analyzing such properties of the proposed variant of SDDP is out of the scope of this paper but will be an interesting future direction to explore.

    \item We numerically compare our approach with LazyPSRL (adapted to the stochastic optimal control problem) and a variant of distributionally robust approach with shrinking ambiguity set. Our approach outperforms LazyPSRL by achieving a smaller regret, and the regret decreases with more samples from the posterior, verifying the intuition of exploiting the entire potesterior is more beneficial in control problems. Our approach also outperforms the distributionally robust approach in the sense that the regret of our approach decreases much faster over time, probably because the distributionally robust approach is too conservative and does not easily adapt its ambiguity set to the observed data.  
\end{itemize}

The rest of the paper is organized as follows. We present the problem formulation and framework of episodic Bayesian SOC in Section 2. In Section 3, we show the asymptotic convergence of episodic value functions and associated policies to the true optimal ones and derive asymptotic rate of episodic value functions. In Section 4, we develop a computational method for episodic Bayesian SOC. In Section 5, we carry out numerical experiments on an inventory control problem to verify our theoretical results and demonstrate the effectiveness of our computational method. Finally, we conclude the paper and outline some future directions in Section 6.

\setcounter{equation}{0}
\section{Episodic Bayesian Optimal Control}
\label{sec-boc}

We consider the discrete-time infinite-horizon  Stochastic Optimal Control (SOC) problem (e.g., \cite{ber78}), as follows
\begin{equation}
\label{exp-1}
\begin{array}{ll}
\min\limits_{\pi\in \Pi}  \bbe^\pi\Big[ \sum_{t=1}^{\infty}\gamma^{t-1}
c(x_t,u_t,\xi_t)
\Big],
\end{array}
\end{equation}
where  variables  $x_t\in \X\subset\bbr^{n}$  represent state  of the system,  $\U$ is a nonempty closed subset of  $\bbr^{m}$,
$u_t\in \U$  are controls,    $c:\X\times \U\times\Xi\to \bbr$ is cost function,  $\gamma\in (0,1)$  is the discount factor,    $F:\X\times\U\times\Xi\to \X$  is a  (measurable) mapping, and $\xi_t$, $t=1,...,$  is an independent identically distributed (iid)  sequence of random vectors viewed   as realizations of random vector $\xi$ whose probability distribution $P_*$ is
supported on the  set $\Xi\subset  \bbr^{d}$ equipped with its Borel sigma algebra $\B$. The optimization (minimization) in \eqref{exp-1} is performed over the set
$\Pi$  of feasible policies
satisfying almost surely (a.s.)
%the constraints $u_t\in \U$  and  the state equations
\begin{equation}\label{exp-2}
u_t\in \U\;\;{\rm and}\;\;x_{t+1}=F(x_t,u_t,\xi_t), \;t\ge 1,
\end{equation}
given initial value $x_1$ of the state vector. It is also possible to consider settings where   the control set $\U(x)$ depends on the state, and we briefly discuss this in the electronic companion \ref{EC: state-dependent control}.

In practice, the true distribution $P_*$ is often unknown and needs to be estimated from data. To continually learn the true distribution from sequential data, it is natural to adopt a Bayesian approach. 
%For computational convenience, 
We assume that the probability distribution of the random vector  $\xi$ is modeled by a parametric family defined  by probability  density function (pdf)  $f(\cdot|\theta)$,
 $\theta\in \Theta$,  where $\Theta\subset \bbr^q$ is a closed set. With the Bayesian perspective, the parameter vector $\theta$  is assumed to be random, whose probability distribution is supported on the set $\Theta$ and defined by a prior probability density $p(\theta)$. Then, given the data (samples) $\bxi^{(N)}=(\xi_1,...,\xi_N)$, the posterior distribution is determined by the Bayes' rule
\begin{equation}\label{stat-2}
 p(\theta|\bxi^{(N)})=\frac{f(\bxi^{(N)}|\theta) p(\theta)}{\int_\Theta f(\bxi^{(N)}|\theta)p(\theta) d\theta},
\end{equation}
where  $f(\bxi^{(N)}|\theta)=\prod_{i=1}^N f(\xi_i|\theta)$ is the  density of the data conditional on $\theta$.

It is possible that our assumed parametric model $f(\cdot|\theta)$ does not include the true distribution $P_*$, but we can still approximate the true distribution by the best approximation within the parametric family. The best approximation is specified by the parameter that minimizes the Kullback--Leibler (KL) divergence from the true distribution to the parametric model (see, e.g., \cite{SZY2023}). More specifically,  assume that the true distribution $P_*$ of $\xi$ has a probability density function (pdf) $q(\cdot)$ and that there is a value $\theta^*\in \Theta$ of the parameter vector which minimizes the KL divergence from $q(\cdot)$ to $f(\cdot|\theta)$. That is,
 \begin{equation}\label{kldiv}
 \theta^*\in \arg\min_{\theta\in \Theta} D_{KL}(q\|f_\theta),
 \end{equation}
where $f_\theta(\cdot):=f(\cdot|\theta)$ and
\[
D_{KL}(q\|f_\theta)=\int_\Xi q(\xi) \ln (q(\xi)/f_\theta(\xi))d\xi.
\]
If $q(\cdot)=f_{\theta^*}(\cdot)$, the parametric model is said to be {\em correctly specified}; otherwise, it is said to be {\em mis-specified}.

The posterior distribution $p(\theta|\bxi^{(N)})$ provides a density estimation of $\theta$ given all the data and can be used to construct a Bayesian average estimate of the unknown true problem.
First, we define the joint distribution  $P_N$ of $(\xi,\theta)$,
% on $(\Xi,\B)$, 
 which is specified by the density $f(\xi|\theta)$ conditional on $\theta$ and the posterior of $\theta$.
Let $\theta_N$ denote a random vector that follows the posterior  $p(\theta|\bxi^{(N)})$.
Thus,  for a random variable $Z:\Xi\to \bbr$,
\begin{equation}\label{post-2a}
\bbe_{P_N}[Z]=\bbe_{\theta_N} \left[\bbe_{|\theta}[Z]\right]=\int_\Theta \int_\Xi Z(\xi)f(\xi|\theta)p(\theta|\bxi^{(N)})d\xi d\theta,
\end{equation}
where
\begin{equation}\label{condexp}
 \bbe_{|\theta}[Z]:=  \int_\Xi Z(\xi)f(\xi|\theta) d\xi
\end{equation}
  is the expectation with respect to distribution  of $\xi$ conditional on $\theta$,
and
\begin{equation}\label{post-2}
\bbe_{\theta_N}[Y]:=\int_\Theta  Y(\theta)p(\theta|\bxi^{(N)})d\theta
\end{equation}
denotes the expectation of random variable $Y:\Theta\to\bbr$  with respect to the posterior distribution $p(\theta|\bxi^{(N)})$. Next, we will use the expectation operator $\bbe_{P_N}[\cdot]$ to define the Bellman equation for the Bayesian average problem in the following.

For the true problem \eqref{exp-1}, the associated Bellman equation for the value function is
 \begin{equation}
 \label{bellm-1}
V(x)=\inf\limits_{u\in \U} \bbe
   \left [c(x,u,\xi)+\gamma
V \big(F(x,u,\xi) \big)\right].
\end{equation}
Assuming    that  the cost function is  bounded,  equation \eqref{bellm-1} has a  unique solution  (e.g., \cite{ber78}).
We denote by $V^*$ the solution of equation \eqref{bellm-1} corresponding to the probability distribution of $\xi$ defined by the pdf $f(\cdot|\theta^*)$, where $\theta^*$ is the parameter vector  which minimizes the  KL-divergence in \eqref{kldiv}. That is
\begin{equation}
 \label{belltheta}
V^*(x)=\inf\limits_{u\in \U} \bbe_{|\theta^*}
   \left [c(x,u,\xi)+\gamma
V^* \big(F(x,u,\xi) \big)\right].
\end{equation}
The  respective   policy is given by $\pi^*: \mathcal{X} \rightarrow \mathcal{U}$ that satisfies  %(e.g., \cite{ber78})
 \begin{equation}\label{bellm-2}
\pi^*(x)\in\arg\min_{u\in \U} \bbe_{|\theta^*}
   \left [c(x,u,\xi)+\gamma
V^* \big(F(x,u,\xi) \big)\right].
\end{equation}
If the model is correctly specified, then $V^*$ is the value function corresponding to the true distribution of $\xi$, and $\pi^*$ is the optimal policy for problem \eqref{exp-1}.

With the expectation operator $\bbe_{P_N}[\cdot]$,  the Bayesian average counterpart of the Bellman equation is
 \begin{equation}
 \label{bellm-3}
V^*_N(x)=\inf\limits_{u\in \U} \bbe_{P_N}
   \big [c(x,u,\xi)+\gamma
V^*_N \big(F(x,u,\xi) \big)\big].
\end{equation}
This dynamic programming equation corresponds to the infinite horizon problem \eqref{exp-1} with the expectation of $\xi$ taken with respect to $P_N$.
For the solution $V^*_N(\cdot)$ of \eqref{bellm-3}, its corresponding policy is determined by
  \begin{equation}
 \label{bellm-4}
\pi_N(x)\in \arg\min\limits_{u\in \U} \bbe_{P_N}
   \big [c(x,u,\xi)+\gamma
V^*_N \big(F(x,u,\xi) \big)\big].
\end{equation}

Since the Bayesian average problem is only an estimate of the true problem given currently available data, we will only exercise the policy $\pi_N(x)$ for one episode. During this episode, we will observe more data of  the randomness $\xi$, and hence, we can use the new data to update the posterior of $\theta$ and the associated policy by resolving the Bellman equation \eqref{bellm-3} and \eqref{bellm-4} with $p(\theta|\bxi^{(N)})$ replaced by the new posterior.
This process generates the exercised policy $\pi = \{\pi_1(\cdot), \pi_2(\cdot), \ldots\}$.

To summarize the process described above, we present the following algorithm of Episodic Bayesian Optimal Control (EBOC) . Suppose that at the start of the first episode, a batch of historic data $\bxi^{(1)} = (\xi_1,...,\xi_{N_0})$ is available.
\\
%{\color{blue}
\begin{minipage}{1\linewidth}
\begin{algorithm}[H]
\SetAlgoLined
\SetKwInOut{Input}{input}\SetKwInOut{Output}{output}
\Input{initial state $x_1$; initial posterior distribution $p(\theta|\bxi^{(1)})$ computed using a historical batch of data $\bxi^{(1)}$.}
%prior distribution $p_0(\theta) = p(\theta|\bxi^{(0)})$.} %upper bound for cost function $\kappa$; discount factor $\gamma$; tolerance for the the error of the truncation $\cT$; initial lower approximation $\underline{V}_1(\cdot)$;
%\Output{policy $\{\pi_1, \pi_2,\ldots\}$}
\For{Episode $N = 1, 2, \ldots$}{
\textbf{Step 1}:
solve the Bellman equation \eqref{bellm-3} and \eqref{bellm-4} to obtain value function $V^*_N(\cdot)$ and corresponding policy $\pi_N$.\\
\textbf{Step 2:} apply $\pi_N$ and observe transition to the next state $x_{N+1} = F(x_N, \pi_N(x_N), \xi_{N+1})$; use the new observation $\xi_{N+1}$ to update the posterior according to
$$
p(\theta|\bxi^{(N+1)}) =  \frac{p(\theta|\bxi^{(N)})f(\xi_N|\theta)}{\int p(\theta|\bxi^{(N+1)})f(\xi_{N+1}|\theta)d\theta},
$$
where $\bxi^{(N+1)} = (\bxi^{N},\xi_{N+1})$.
}
\caption{Episodic Bayesian optimal control (EBOC)}
\label{algorithm: BayesianOC}
\end{algorithm}
\end{minipage}
%}
Please note that here we implicitly assume each episode has length $1$, i.e., each current policy is only exercised for one step. However, our approach can be easily generalized to any deterministic positive integer $H$, where each current policy is used for $H$ steps. During such an episode, $H$ data points are collected and then used to update the posterior and re-solve the Bellman equation under the new posterior. Algorithm~\ref{algorithm: BayesianOC} is mostly conceptual, since its Step 1 is usually not exactly solvable except for some special and simple cases. We will investigate computational methods for Step 1 (i.e., solving the Bayesian average Bellman equation) in Section~\ref{sec-comp}.
%}

%\setcounter{equation}{0}
\section{Convergence Analysis} \label{sec:convergence}

Our proposed Episodic Bayesian Optimal Control approach generates a sequence of episodic value functions and policies. In this section, we will analyze their convergence with respect to $N$, the number of episodes or data size.

\subsection{Convergence of the episodic value functions and policies}
%Let $P_*$ be the probability distribution supported on  $(\Xi,\B)$ defined by the pdf  $f(\cdot|\theta^*)$ with $\theta^*\in \Theta$ being the true value of the parameter vector. That is,    for a random variable $Z:\Xi\to \bbr$ we have $\bbe_{P_* }[Z]=\int_\Xi Z(\xi) f(\xi|\theta^*)d\xi$.

Denote by
$\bbb(\X)$ the space of bounded      functions $v:\X\to\bbr$ equipped with the sup-norm $\|v\|_\infty:=\sup_{x\in \X}|v(x)|$.
Unless stated otherwise, probabilistic statements like almost surely (a.s.) are made with respect to the true distribution  $P_*$ of $\xi$. We make the following assumptions.

 \begin{assumption}
 \label{ass-1}
{\rm (i)} The cost function
$c(x,u,\xi)$ is bounded on  $\X\times \U\times \Xi$. {\rm (ii)}  For every $\xi\in \Xi$ the function $f(\xi|\cdot)$ is continuous at $\theta^*$, and
there is a (measurable) function $h:\Xi\to \bbr_+$ such that $\int_\Xi h(\xi)d\xi$ is finite and $f(\xi|\theta)\le h(\xi)$ for all $\xi\in \Xi$ and all $\theta$ in a neighborhood of $\theta^*$.
\end{assumption}

 \begin{assumption}{\bf(Bayesian Consistency)}
 \label{ass-2} Almost surely
   $\theta_N$ converges in probability to $\theta^*$, i.e., for almost every sequence  $\{\xi_1,...\}$ and  any $\e>0$   it follows that
 \begin{equation}\label{LLN}
 \lim_{N\to\infty} \int_{ \theta:\|\theta-\theta^*\|\le \e}   p(\theta|\bxi^{(N)}) d\theta = 1.
 \end{equation}
 \end{assumption}
Intuitively, Assumption~\ref{ass-2} means that the Bayesian posterior of $\theta$ becomes more concentrated on $\theta^*$ as the data size increases and eventually degenerates to the delta function on $\theta^*$. There are classical results, going back to   \cite{doob1948application} and  \cite{Schwartz1965OnBP},  ensuring Bayesian consistency (i.e., Assumption \ref{ass-2}) under certain regularity conditions.  Relatively simple regularity conditions for such results were suggested in \cite{SZY2023}.   For  reader's convenience these regularity conditions are given below.

\paragraph{Regularity conditions for Assumption~\ref{ass-2}.}\label{EC:assump2}

%\begin{assumption} \label{assump}
{\rm (i)} The set $\Theta$ is convex  compact with nonempty interior.
{\rm (ii)} $\ln p(\theta)$ is bounded on $\Theta$, i.e., there are constants $c_1>c_2>0$ such that $c_1\ge p(\theta)\ge c_2$ for all $\theta\in \Theta$.
{\rm (iii)}  $q(\xi)>0$  for $\xi\in \Xi$.
{\rm (iv)} $f(\xi|\theta)>0$, and hence  $p(\theta|\bxi^{(N)})>0$,  for all $\xi\in \Xi$ and $\theta\in \Theta$. {\rm (v)}     $f(\xi|\theta)$ is continuous in $\theta\in \Theta$.
{\rm (vi)} $\ln f(\xi|\theta)$, $\theta\in \Theta$,  is dominated by an integrable (with respect to the true distribution)  function.
%\end{assumption}

%Let $V^*$ be the solution of equation \eqref{bellm-1} where the expectation is taken with respect to $P_*$.
 Assumption \ref{ass-1}(i) guarantees the existence and uniqueness of solutions $V^*$ and $V^*_N$ of the respective  Bellman equations. It also follows that  $V^*$ and $V^*_N$ are bounded, i.e., $V^*,V^*_N\in \bbb(\X)$.

 \begin{proposition}
 \label{pr-1}
 Suppose that Assumptions {\rm \ref{ass-1}}  and  {\rm \ref{ass-2}} hold. Then,  a.s. solution $V^*_N$ of equation \eqref{bellm-3}
converges  uniformly  to $V^*$,   i.e.,
\begin{equation}\label{limit}
 \|V^*_N-V^*\|_\infty\to 0\;\;a.s.,~ \text{as}~ N\to \infty.
\end{equation}
 \end{proposition}

{\bf Proof.}
%Let   $P_N$ be the probability distribution on $\Xi$ corresponding to the posterior distribution,  and hence  $\bbe_{P_N}[Z]=\bbe_{\theta_N}\left [\bbe_{|\theta}[Z]\right]$.
%\begin{equation}\label{pn}
% \bbe_{P_N}[X]=\bbe_{\theta_N}\left [\bbe_{|\theta}[X]\right]=\int_\Theta \left(\int_\Xi X(\xi) f(\xi|\theta)d\xi)\right)p(\theta|\bxi^{(N)})d\theta.
%\end{equation}
The convergences below are understood  in the a.s. sense,
i.e., for almost every (a.e.) $\{\xi_1,...,\}$.
 Consider Bellman operators $\cT,\cT_N:\bbb(\X)\to\bbb(\X)$  defined as
\begin{eqnarray}
\label{operat-1}
\cT(V)(x)&:=&\inf\limits_{u\in \U} \bbe_{|\theta^*}
   \left [c(x,u,\xi)+\gamma
V \big(F(x,u,\xi) \big)\right],\;\;V\in \bbb(\X),\\
\cT_N(V)(x)&:=&\inf\limits_{u\in \U} \bbe_{P_N}
   \left [c(x,u,\xi)+\gamma
V \big(F(x,u,\xi) \big)\right], \;\;V\in \bbb(\X).
\label{operat-2}
\end{eqnarray}
Note that  $\cT$ and $\cT_N$ are contraction mappings (e.g.,    \cite{ber78}),  i.e.,
\begin{equation}\label{contr}
\|\cT (V)-\cT(V')\|_\infty\le \gamma \|V-V'\|_\infty,\;\;V,V'\in \bbb(\X),
\end{equation}
and similarly for $\cT_N$.  Hence, $V^*$
 is the fixed point of $\cT$, i.e., $\cT(V^*)=V^*$, and  $V^*_N$ is the fixed point of $\cT_N$.
 We have
 \begin{equation}\label{valineq}
  \begin{array}{lll}
\|V^*_N-V^*\|_\infty&=&\|\cT_N(V^*_N)-\cT(V^*)\|_\infty\\
&\le& \|\cT_N(V^*_N)-\cT_N(V^*)\|_\infty+\|\cT_N(V^*)-\cT(V^*)\|_\infty\\
&\le & \gamma \|V^*_N-V^*\|_\infty+\|\cT_N(V^*)-\cT(V^*)\|_\infty,
\end{array}
 \end{equation}
 where the last inequality follows by the contraction property of $\cT_N$.
 Thus
 \begin{equation}
 \label{eq-est}
 \|V^*_N-V^*\|_\infty\le (1-\gamma)^{-1} \|\cT_N(V^*)-\cT(V^*)\|_\infty.
 \end{equation}
 It follows from \eqref{eq-est} that  in order to show that $V^*_N(x)$ converges to $V^*(x)$ uniformly in $x\in \X$,   it suffices to show that   $\cT_N(V^*)$ converges to $\cT(V^*)$ with respect to the sup-norm.    Consequently  by \eqref{operat-1} and \eqref{operat-2}
 it suffices  to show that
 $\inf_{u\in \U}\bbe_{P_N} [Z_{(x,u)}]$ converges uniformly to
$\inf_{u\in \U}\bbe_{|\theta^*} [Z_{(x,u)}]$ in $x\in \X$, where
\begin{equation}\label{zfun}
Z_{(x,u)}(\xi):=c(x,u,\xi)+\gamma V^*\big(F(x,u,\xi) \big).
\end{equation}
  In turn, for that it suffices to show that a.s.
 $\bbe_{P_N}  \left [Z_{(x,u)}\right]$ converges to
 $\bbe_{|\theta^*}\left[Z_{(x,u)}\right]$
 uniformly in $(x,u)\in \X\times \U$   (see Lemma \ref{lem-uniform} in the electronic companion \ref{sec:EClemma}).

Since $c(\cdot,\cdot,\cdot)$ and $V^*(\cdot)$  are  bounded,   there is $C>0$   such  that $|Z_{(x,u)}(\xi)|\le C$ for all $(x,u,\xi)\in \X\times \U\times \Xi$. Then    by     \eqref{post-2a} and since
\begin{equation}\label{thetastar}
\bbe_{|\theta^*}[Z_{(x,u)}]=  \int_\Xi  \int_{\Theta} Z_{(x,u)}(\xi) f(\xi|\theta^*) p(\theta|\bxi^{(N)})d\theta d\xi,
\end{equation}
we have  that
   \begin{eqnarray}
  \nonumber  \big | \bbe_{P_N}[Z_{(x,u)}]- \bbe_{|\theta^*}[Z_{(x,u)}]\big|&=&
    \left| \int_{\Theta}  \int_\Xi Z_{(x,u)}(\xi)\big ( f(\xi|\theta)-f(\xi|\theta^*) \big) p(\theta|\bxi^{(N)})d\xi d\theta\right|\\
  \nonumber   & \le & C   \int_{\Theta}   \int_\Xi   \big | f(\xi|\theta)-f(\xi|\theta^*) \big | p(\theta|\bxi^{(N)}) d\xi  d\theta\\
    &=& C   \int_{\Theta}  \psi(\theta) p(\theta|\bxi^{(N)})  d\theta
    =C(\nu_1(\e) +\nu_2(\e)),
    \label{exp-1a}
   \end{eqnarray}
   where
   $
 \psi(\theta):=\int_\Xi  \big| f(\xi|\theta)-f(\xi|\theta^*) \big|  d\xi,
 $
and for $\e>0$,
\begin{equation}
\label{funv}
\nu_1(\e):=\int_{\theta:\|\theta-\theta^*\|\le \e}  \psi(\theta) p(\theta|\bxi^{(N)})  d\theta,\;\;
\nu_2(\e):=\int_{\theta:\|\theta-\theta^*\|> \e}  \psi(\theta) p(\theta|\bxi^{(N)})  d\theta.
\end{equation}

We have that $0\le \psi(\theta)\le  \int_\Xi    f(\xi|\theta) d\xi+\int_\Xi    f(\xi|\theta^*) d\xi=2$, that $\psi(\theta^*)=0$, and by Assumption \ref{ass-2}(ii)
\[
\lim_{\theta\to \theta^*}\psi(\theta)=\int_\Xi \lim_{\theta\to \theta^*} \big| f(\xi|\theta)-f(\xi|\theta^*) \big|  d\xi=0,
\]
where
the interchange of the limit and integration follows by
 Lebesgue's dominated convergence theorem.
Therefore, $\sup\{\psi(\theta): \|\theta-\theta^*\|\le \e\}$ can be made arbitrarily small for sufficiently small $\e>0$, and hence the term $\nu_1(\e)$ can be made arbitrarily small. The term $\nu_2(\e)$ can be bounded  as
$$\nu_2(\e)\le 2 \int_{\theta:\|\theta-\theta^*\|> \e}   p(\theta|\bxi^{(N)})  d\theta,$$   and hence by \eqref{LLN},  a.s. for any  $\e>0$  this term tends to zero as $N\to \infty$. Therefore, the conclusion \eqref{limit} follows.
$\hfill\square$
\\

Uniform convergence of $V^*_N$ to $V^*$ implies convergence of policy $\pi_N$, defined by \eqref{bellm-4}, to a policy    of the limiting problem, defined by  \eqref{bellm-2}. That is, consider the set
\begin{equation}
 \U^*(x):=\arg\min_{u\in \U} \bbe_{|\theta^*}
   \left [c(x,u,\xi)+\gamma
V^* \big(F(x,u,\xi) \big)\right],\;\;x\in \X.
\end{equation}
If the model is correctly specified, then
  $\pi^*(x)\in \U^*(x)$ defines an optimal policy for the  true   SOC problem \eqref{exp-1}.

\begin{proposition}
\label{pr-polcon}
Suppose that Assumptions {\rm \ref{ass-1}}  and  {\rm \ref{ass-2}} hold, and that for  a given $\bar{x}\in \X$ the set $\U^*(\bar{x})$ is nonempty and bounded.  Then, a.s. the distance from $\pi_N(\bar{x})$ to the set  $\U^*(\bar{x})$ tends to zero
as $N\to\infty$. In particular, if $\U^*(\bar{x})=\{\pi^*(\bar{x})\}$  is the singleton, then
a.s. $\pi_N(\bar{x})$ converges to $\pi^*(\bar{x})$.
\end{proposition}

{\bf Proof.}
Consider
\begin{eqnarray*}
\varphi_N(x,u)&:=&\bbe_{P_N}\big [c(x,u,\xi)+\gamma
V^*_N \big(F(x,u,\xi) \big)\big]\\
\varphi(x,u)&:=&\bbe_{|\theta^*}\big [c(x,u,\xi)+\gamma
V^* \big(F(x,u,\xi) \big)\big].
\end{eqnarray*}
 Consider an element $\pi(\bar{x})$ of the set  $\U^*(\bar{x})$. Note that $\pi(\bar{x})\in \arg\min_{u\in \U} \varphi(\bar{x},u)$ and
 \begin{equation}\label{eqpi}
  \pi_N(\bar{x})\in \arg\min_{u\in \U} \varphi_N(\bar{x},u).
 \end{equation}
 By \eqref{limit} and the proof of Proposition \ref{pr-1}  (compare with \eqref{exp-1a} and \eqref{funv}) we have that
\[
\sup_{(x,u)\in \X\times \U} |\varphi_N(x,u)-\varphi(x,u)|\to 0\;\;a.s.,~ \text{as}~ N\to \infty.
\]
Therefore, a.s. for any $\e>0$ there is   $N_\e$ such that for any $N\ge N_\e$ it follows that
\begin{eqnarray*}
\varphi(\bar{x},\pi_N(x))-\e/2\le \varphi_N(\bar{x},\pi_N(x))\\
\varphi_N(\bar{x},\pi(x)) \le \varphi(\bar{x},\pi(x))+\e/2.
\end{eqnarray*}
Moreover, since $\pi(\bar{x})\in \U$ and because of \eqref{eqpi} we have that
\[
\varphi_N(\bar{x},\pi_N(\bar{x}))\le \varphi_N(\bar{x},\pi(\bar{x})).
\]
That is,
\[
 \varphi(\bar{x},\pi_N(\bar{x}))\le \varphi(\bar{x},\pi(\bar{x}))+\e,
 \]
i.e.,    $\pi_N(\bar{x}) \in L_{\e}(\bar{x})$, where $L_\e(\bar{x})$ is the level set
\[
  L_\e(\bar{x}):=\{u\in \U: \varphi(\bar{x},u)\le  \varphi(\bar{x},\pi(\bar{x}))+\e\}.
\]
Now let $\e_k$ be a decreasing sequence of positive numbers  converging to zero. Note that $L_{\e_{k+1}}(\bar{x})\subseteq L_{\e_k}(\bar{x})$ and
$\cap_{k=1}^\infty L_{\e_k}(\bar{x})$ is contained in the topological closure of the set $\U^*(\bar{x})$. Indeed,   otherwise there is a sequence $u_k\in L_{\e_k}(\bar{x})$ such that the distance   $\dist(u_k,\U^*(\bar{x}))\ge \delta$  for some $\delta>0$. Since $\U^*(\bar{x})$ is bounded we can choose such sequence to be bounded, and hence it has an accumulation point $u^*$. It follows that  $\varphi(\bar{x},u^*)\le \varphi(\bar{x},\pi(\bar{x}))+\e$ for any $\e>0$, and hence  $\varphi(\bar{x},u^*)= \varphi(\bar{x},\pi(\bar{x}))$, i.e. $u^*\in \U^*(\bar{x})$. On the other hand,
$\dist(u^*,\U^*(\bar{x}))\ge \delta$, a contradiction.

It follows that a.s. the distance from $\pi_N(\bar{x})$ to the topological closure of $\U^*(\bar{x})$ tends to zero, and hence the distance from $\pi_N(\bar{x})$ to  $\U^*(\bar{x})$ tends to zero.
$\hfill\square$
\\

We see that under mild  regularity  conditions, the sequence of value functions $V^*_N$, and the corresponding policies $\pi_N$, determined by   Bellman equation \eqref{bellm-3}, converge  to their counterparts of the limit  problem. If the parametric model is correctly specified for the true distribution, the limit problem corresponds to the true problem.
The next natural question is how fast is such convergence. For this we will derive the asymptotics of the episodic value functions $V^*_N$ in the following subsection.

%Our derivation is somewhat heuristic, while a rigorous derivation of such results is far beyond the scope of this paper. These heuristics  are verified in the numerical experiments of section 5.
%Our derivations are  based on calculation of the Influence Function (see the approximation \eqref{appr-1} below). Such  {\color{blue}   approximations} are  routinely used in statistical applications (cf., \cite{Fernholz},\cite[Section 20.1]{vaart}).
%, while a rigorous derivation of such results is far beyond the scope of this paper.
%Further the derived asymptotics   are verified in the numerical experiments of section 5.

\subsection{Asymptotics  of  episodic value functions} \label{sec-asymptotic}

%Our derivations are  based on calculation of the Influence Function (see the approximation \eqref{appr-1} below). Such  {\color{blue}   approximations} are  routinely used in statistical applications (cf., \cite{Fernholz},\cite[Section 20.1]{vaart}). Further the derived asymptotics are verified in the numerical experiments of section 5.

In this subsection, we want to study the asymptotic convergence rate of the episodic value function, i.e, the asymptotics of the value function $V^*_N(x)$ for a given fixed value $x=x'$.  We use the approach of (von Mises) statistical functionals. Specifically,  denote by $\vv(P)$ the optimal value of problem \eqref{exp-1}, considered as a function of the distribution $P$ of $\xi_t$ for initial value $x_1=x'$. We have that $\vv(P_*)=V^*(x')$ and $\vv(P_N)=V^*_N(x')$. Consider the directional derivative
\begin{equation} \label{dirdifvalue}
  \vv'(P_*,P-P_*)=\lim_{\tau\downarrow 0}\frac{\vv(P_*+\tau(P-P_*))-\vv(P_*)}{\tau}.
\end{equation}
If  $\vv'(P_*,\cdot)$ exists and moreover  is linear, then $\vv(\cdot)$ is {\em G\^ateaux }  differentiable  at $P_*$. In that case  $\vv'(P_*,P_N-P_*)$ is  called the {\em  Influence Function} of $\vv(\cdot)$.

The directional derivative \eqref{dirdifvalue} is evaluated as follows. Let $\pi^*(x)$ be an optimal solution of  problem \eqref{exp-1} for $P=P^*$. Then $\vv(P_*)=\bbe_{P_*}  \Big[ \sum_{t=1}^{\infty}\gamma^{t-1}
c(x_t,\pi^*(x_t),\xi_t)
\Big]$. Since the policy $\pi^*(x)$ is feasible, we have that 
\[
\vv(P_*+\tau(P-P_*))\le \bbe_{P_*+\tau(P-P_*)}  \Big[ \sum_{t=1}^{\infty}\gamma^{t-1}
c(x_t,\pi^*(x_t),\xi_t)
\Big].
\]
By linearity of the expectation operator, it follows for $\tau\in (0,1)$ that 
\begin{equation}\label{dirdifvalue-2}
  \frac{\vv(P_*+\tau(P-P_*))-\vv(P_*)}{\tau}\le
 \bbe_{P-P_*}  \Big[ \sum_{t=1}^{\infty}\gamma^{t-1}
c(x_t,\pi^*(x_t),\xi_t)
\Big] .
\end{equation}
To derive the  lower bound, let $\pi_\tau(x)$ be an optimal solution of \eqref{exp-1} for $P_*+\tau(P-P_*)$. Then by similar arguments, we have 
\begin{equation}\label{dirdifvalue-3}
 \frac{\vv(P_*+\tau(P-P_*))-\vv(P_*)}{\tau}\ge
 \bbe_{P-P_*}  \Big[ \sum_{t=1}^{\infty}\gamma^{t-1}
c(x_t,\pi_\tau(x_t),\xi_t)
\Big] .
\end{equation}
With the two inequalities above, we can conclude that 
\begin{equation}\label{direc}
\vv'(P_*,P-P_*)=\bbe_{P-P_*}  \Big[ \sum_{t=1}^{\infty}\gamma^{t-1}
c(x_t,\pi^*(x_t),\xi_t)
\Big],
\end{equation}
provided  the right hand side of \eqref{dirdifvalue-3} converges to the right hand side of \eqref{dirdifvalue-2} as $\tau\downarrow 0$. %, which is the main technical difficulty for establishing \eqref{direc}. 
Proving such convergence is the main technical difficulty of verifying \eqref{direc}. For that we certainly need to assume that for $P=P_*$, problem \eqref{exp-1} has a unique optimal solution $\pi^*$, i.e., $\pi^*(x)$ is the unique minimizer in the right hand side of \eqref{bellm-2} for every  $x\in \X$. Under the uniqueness assumption, the convergence can be proved as in the Danskin Theorem (see Section \ref{danskin} in the Electronic Companion), which additionally assumes compactness of the feasible set. Such compactness can be verified in some cases,  in particular when either the set $\X$ is finite and the set $\U$  is compact, or the set $\U$ is finite. More generally, the required convergence can be verified on a case-by-case basis.  For instance, in our numerical example of an inventory control problem with uncountable sets $\X$ and $\U$, the asymptotics result (Proposition 3) holds well (see Section~\ref{sec:numerical}).

% Suppose that for $P=P_*$  problem \eqref{exp-1} has unique optimal solution $\pi^*$, i.e., $\pi^*(x)$ is the unique minimizer in the right hand side of \eqref{bellm-2} for every  $x\in \X$. Denote  $\bbe_{P-P_*}[Z] :=\bbe_{P}[Z]-\bbe_{P_*}[Z]$. Then under certain regularity conditions,
% \begin{equation}\label{direc}
% \vv'(P_*,P-P_*)=\bbe_{P-P_*}  \Big[ \sum_{t=1}^{\infty}\gamma^{t-1}
% c(x_t,\pi^*(x_t),\xi_t)
% \Big],
% \end{equation}
% with $x_1=x'$ and $x_{t+1}=F(x_t,\pi^*(x_t),\xi_t)$, $t\ge 1$.  That is, for $P=P_N$ the right hand side of \eqref{direc} is the Influence Function of $\vv(\cdot)$.
% Formula \eqref{direc} is a counterpart of such formula in the setting of  static stochastic programming problems (see \cite[Theorem 5.7 and eq. (5.27)]{SDR}), which is based on
% the Danskin Theorem (see  the details in electronic companion \ref{danskin}). The main technical difficulty in applying such results is verification of compactness of the respective feasible set. This is
% discussed in \cite{shacheng2021}.  In particular \eqref{direc} holds if either  the set $\X$ is finite and the set $\U$  is compact, or the set $\U$ is finite. {\color{blue} However, these assumptions can be relaxed in practice; in particular, we numerically verified our theoretical results on an inventory control problem with uncountable sets $\X$ and $\U$, and the results hold well (see Section~\ref{sec:numerical}).}

Now we use the von Mises expansion (e.g., Section 20.1 in \cite{vaart}), which can be intuitively understood as the first-order approximation of the functional:
\begin{equation}\label{appr-1}
  \vv(P_N)-\vv(P_*)\approx \vv'(P_*,P_N-P_*).
\end{equation}
For conditions under which this expansion is meaningful in terms of the convergence of remainders, we refer interested readers to \cite{Fernholz1983vonMises} for a detailed discussion.

This, together with \eqref{direc}, gives us  the following approximation
\begin{equation}\label{appr-2}
\begin{array}{lll}
 V^*_N(x')-V^*(x')&\approx& \bbe_{P_N-P_*} \Big[ \sum_{t=1}^{\infty}\gamma^{t-1}
c(x_t,\pi^*(x_t),\xi_t)\Big]\\
&=&\bbe_{P_N} \Big[ \sum_{t=1}^{\infty}\gamma^{t-1}
c(x_t,\pi^*(x_t),\xi_t)\Big]-V^*(x').
\end{array}
\end{equation}
Since we only care about the convergence rate of  $ V^*_N(x')$ to $V^*(x')$ and the rate is dominated by the first-order term, this approximation will not affect the convergence rate.

By \eqref{appr-2}, it suffices to consider the asymptotic convergence rate of 
\begin{equation}\label{PN}
\bbe_{P_N}[Z]- \bbe_{P_*}[Z]=
    \int_\Xi  \int_{\Theta} Z (\xi) \big ( f(\xi|\theta)-f(\xi|\theta^*) \big) p(\theta|\bxi^{(N)}) d\theta d\xi.
\end{equation}
%where $Z$ is set to $\sum_{t=1}^{\infty}\gamma^{t-1} c(x_t,\pi^*(x_t),\xi_t)$. 
Note that $\bbe_{P_N}[Z]- \bbe_{P_*}[Z]$ is a function of the random vector $\bxi^{(N)}$ and its convergence rate depends on the convergence rate of the Bayesian posterior $P_N$ to the true probability distribution of $\xi$, $P_*$. Therefore, we will use the Bernstein-von Mises theorem, which is also known as the Bayesian central limit theorem. Specifically, we assume in this section that the model is {\em correctly specified}; i.e., $P_*$ is defined by the pdf $f(\cdot|\theta^*)$, and hence $\bbe_{P_*}=\bbe_{|\theta^*}$.
The Bernstein-von Mises theorem states that under certain regularity conditions, in particular  that the prior is continuous  and (strictly) positive in a neighborhood around $\theta^*$, then
\begin{equation}\label{bern}
\left\|p(\theta|\bxi^{(N)})-
d\hspace{0.5mm}
\N\big(\hat{\theta}_N,N^{-1}I(\theta^*)^{-1}\big )\right\|_{TV}\stackrel{P_{\theta^*}^N}{\to} 0.
\end{equation}
Here $p(\theta|\bxi^{(N)})$ is the pdf of the posterior distribution $P_N$,
$ I(\theta)=\bbe_{P_*}\left[S(\theta) S(\theta)^\top\right ]$
  is the Fisher information matrix,
  \begin{equation}\label{score}
 S(\theta)(\xi):=\nabla_\theta \log f(\xi|\theta)
\end{equation}
 is the score function,
  $\hat{\theta}_N$ is an asymptotically efficient estimator of $\theta^*$  (i.e.,  $N^{1/2}(\hat{\theta}_N-\theta^*)$  converges in distribution to normal $\N(0,I(\theta^*)^{-1})$),
$d\hspace{0.5mm}\N(\mu,\Sigma)$ denotes the pdf of the normal distribution with mean $\mu$ and covariance $\Sigma$,   and
\begin{equation}\label{TV}
\|h-g\|_{TV}=\half \int |h(\theta)-g(\theta)|d\theta
\end{equation}
  is the  total-variation (TV) distance between two pdfs $h$ and $g$.  For a proof of the  Bernstein - von Mises theorem and discussion of the involved regularity conditions we can refer to   \cite[Theorem 10.1 and discussion on page  144]{vaart}).
%In the proposition \ref{pr-convergence} below we assume that the  Bernstein - von Mises theorem is valid, i.e., the convergence \eqref{bern} holds.}

%This, {\color{blue} together with the asymptotic result of Proposition \ref{pr-convergence}},  suggests that the stochastic error of  $V^*_N(x')$   considered as an estimator of $V^*(x')$ is of order $O_p(N^{-1/2})$.

% may need this later:  Recall that a sequence $X_n$ of random variables converges in {\em total variation} to a random variable $X$ if
% \[
% \lim_{n\to\infty} \sup_{A}|P(X_n\in A)- P(X\in A)|=0,
% \]
% where the supremum is taken over all measurable sets $A$. If $X_n$ and $X$ have respective densities $g_n$ and $g$  with respect to $P$, then
% \[
% \sup_{A}|P(X_n\in A)- P(X\in A)|=\half \int |g_n-g|dP,
% \]
% (e.g.,  \cite[Page 22]{vaart}).  Note that convergence in total variation implies convergence in distribution.

To show the limit of \eqref{PN}, we also need to assume the limit and  integration in the right hand side of  \eqref{PN}  can be interchanged, i.e., 
\begin{equation}\label{limitinter}
\lim_{N\to\infty}\left(\bbe_{P_N}[Z]- \bbe_{P_*}[Z]\right)=
\int_\Xi  \int_{\Theta}\lim_{N\to\infty}  W_N (\xi,\theta)   d \theta  d\xi,
\end{equation}
where $W_N(\xi,\theta):= Z (\xi) \big ( f(\xi|\theta)-f(\xi|\theta^*) \big)p(\theta|\bxi^{(N)})$, and the limit is understood almost surely (a.s.) with respect to the true distribution of the data, $P_*$. Conditions for this interchangeability can be found in the electronic companion \ref{sec:ECinterchange}.

The following proposition formally states the asymptotics of \eqref{PN}. 

 \begin{proposition}
 \label{pr-convergence}
 Suppose that (i) the convergence \eqref{bern} holds (Bernstein-von Mises theorem), (ii) $f(\xi|\cdot)$ is differentiable at $\theta^*$ for every $\xi\in \Xi$, (iii) $\theta^*$ is an interior point of $\Theta$, and (iv) the
 limit  and integration  can be interchanged (i.e., \eqref{limitinter} holds a.s.).  Then,   almost surely
   $N^{1/2}\big(\bbe_{P_N}[Z]- \bbe_{P_*}[Z]\big)$ converges in  distribution    to $\N(0,\sigma^2)$
   with
\begin{equation}
\label{sigm}
   \sigma^2= \var\left(  X^\top  \bbe_{P_*}[Z  S(\theta^*)]
 \right)
 =
  \bbe_{P_*}[Z S(\theta^*)]^\top I(\theta^*)^{-1}
 \bbe_{P_*}[Z  S(\theta^*)],
\end{equation}
 where $X\sim \N(0,I(\theta^*)^{-1})$ is a normally distributed random vector and $S(\theta)$ is the score function.
\end{proposition}

{\bf Proof.}
Denote by $\phi(\cdot)$   the pdf of $\N(0,I(\theta^*)^{-1})$.
Then the pdf of  normal distribution with mean $\hat{\theta}_N$ and covariance matrix $N^{-1}I(\theta^*)^{-1}$ is $N^{1/2} \phi\big(N^{1/2}(\theta-\hat{\theta}_N)\big)$, i.e.,  the pdf $d\hspace{0.5mm}
\N\big(\hat{\theta}_N,N^{-1}I(\theta^*)^{-1})$ in \eqref{bern} can be expressed as $N^{1/2} \phi\big(N^{1/2}(\theta-\hat{\theta}_N)\big)$. For simplicity, denote this pdf by $\varphi_N(\theta)$. 

Let's rewrite
\begin{eqnarray} \label{eq-limit}
  && \bbe_{P_N}[Z]- \bbe_{P_*}[Z]  =
   \int_\Xi  \int_{\Theta} Z (\xi) \big ( f(\xi|\theta)-f(\xi|\theta^*) \big) p(\theta|\bxi^{(N)}) d\theta d\xi = \nonumber \\
   && \!\!\!\! \underbrace{\int_\Xi  \int_{\Theta} Z (\xi) \big ( f(\xi|\theta)-f(\xi|\theta^*) \big)\big(p(\theta|\bxi^{(N)}) - \varphi_N(\theta) \big) d\theta d\xi}_{= (A)} 
+ \underbrace{\int_\Xi  \int_{\Theta} Z (\xi) \big ( f(\xi|\theta)-f(\xi|\theta^*) \big) \varphi_N(\theta) d\theta d\xi}_{=(B)}. 
\end{eqnarray}
Taking the limit $N\rightarrow \infty$ on \eqref{eq-limit}, we note that the first term $(A)$ vanishes because
\begin{eqnarray*}
&&\left|\lim_{N\rightarrow\infty}\int_\Xi  \int_{\Theta} Z (\xi) \big ( f(\xi|\theta)-f(\xi|\theta^*) \big)\big(p(\theta|\bxi^{(N)}) - \varphi_N(\theta) \big) d\theta d\xi \right| \\
&=& \left|\int_\Xi Z (\xi) \lim_{N\rightarrow\infty} \int_{\Theta} \big ( f(\xi|\theta)-f(\xi|\theta^*) \big)\big(p(\theta|\bxi^{(N)}) - \varphi_N(\theta) \big) d\theta d\xi \right| \\
&\leq& 2 \sup_{\Theta}\big| f(\xi|\theta)-f(\xi|\theta^*) \big| \left(\lim_{N\rightarrow\infty} \|p(\theta|\bxi^{(N)}) - \varphi_N(\theta) \|_{TV}\right) \int_\Xi | Z (\xi) | d\xi  \\
&=& 0,
\end{eqnarray*}
where the first equation follows from the interchangeability of limit and integration, and the last equation uses the Bernstein-von Mises theorem \eqref{bern}. 

So we only need to focus on the second term $(B)$ of \eqref{eq-limit}. By making the change of variable  $\theta'=N^{1/2}(\theta-\theta^*)$ and defining 
\[
G_N(\xi,\theta):= \frac{ f(\xi|N^{-1/2}\theta+\theta^*)-f(\xi|\theta^*)}{N^{-1/2}}, 
\]
we can write $(B)$ as
\begin{eqnarray}
N^{-1/2} \int_\Xi \int_{\Theta'}
Z (\xi)G_N(\xi,\theta) \phi\big(\theta-N^{1/2}(\hat{\theta}_N-\theta^*))d\theta d\xi,
%= \int_\Xi \int_{\Theta'}Z (\xi)G_N(\xi,\theta) \phi\big(\theta-X)d\theta d\xi,
\label{approx-1}
\end{eqnarray}
where $\Theta':=N^{1/2}(\Theta-\theta^*)$.

Next, observe that $G_N(\xi,\theta)$ tends to $\theta^\top \nabla_\theta f(\xi|\theta^*)$ as $N\to \infty$, provided $f(\xi|\cdot)$ is differentiable at $\theta^*$; and that $N^{1/2}(\hat{\theta}_N-\theta^*)$ converges in distribution to $\N(0,I(\theta^*)^{-1})$, which is the distribution of the random vector $X$.  Recall  that it is assumed that $\theta^*$ is an interior point of the set $\Theta$, and hence  the set $\Theta'$  approaches  $\bbr^q$ as $N\to \infty$ (in particular if $\Theta=\bbr^q$, then  $\Theta'=\bbr^q$ as well). Then, by the interchangeability of integration and limit, we have
\begin{eqnarray*}
\lim_{N\rightarrow\infty}  \int_\Xi \int_{\bbr^q} Z (\xi)G_N(\xi,\theta) \phi\big(\theta-N^{1/2}(\hat{\theta}_N-\theta^*)\big)d\theta d\xi 
= \int_\Xi \int_{\bbr^q} Z (\xi)\theta^\top \nabla_\theta f(\xi|\theta^*) \phi(\theta-X)d\theta d\xi.
\end{eqnarray*}

By making another change of variable $\theta' = \theta-X$ and noticing that $
  \int_{\bbr^q}
 \theta^\top \nabla_\theta f(\xi|\theta^*) \phi\big(\theta)d\theta=0
$, we can write 
\begin{eqnarray}
\nonumber
 &&\int_\Xi \int_{\bbr^q}
Z (\xi)\theta^\top \nabla_\theta f(\xi|\theta^*)\phi\big(\theta-X)d\theta d\xi
  \\
  &=& \int_\Xi \int_{\bbr^q}
Z (\xi)(\theta+X)^\top \nabla_\theta f(\xi|\theta^*)\phi\big(\theta)d\theta d\xi \nonumber\\
 \nonumber  &=& X^\top\left(\int_\Xi
Z (\xi)\big[ \nabla_\theta f(\xi|\theta^*)\big]  d\xi\right)
  \nonumber     \\
  &=&X^\top\left(\int_\Xi
Z (\xi) S(\theta^*)(\xi)f(\xi|\theta^*)\big]  d\xi\right) \nonumber\\
&=& X^\top \bbe_{P_*}[Z S(\theta^*)].
\label{zapprox}
\end{eqnarray}
Therefore, $N^{1/2}\left(\bbe_{P_N}[Z]- \bbe_{P_*}[Z]\right)$ converges to $X^\top\bbe_{P_*}[Z S(\theta^*)]$ a.s. as $N\rightarrow\infty$, and this completes the proof.    $\hfill\square$

Setting $Z =  \sum_{t=1}^{\infty}\gamma^{t-1}
c(x_t,\pi^*(x_t),\xi_t)$, Proposition~\ref{pr-convergence}
implies that the stochastic error of $V_N^*(x')$, 
considered as an estimator of $V^*(x')$, is of the order $O_p(N^{-1/2})$.

Proposition~\ref{pr-convergence} not only shows the asymptotic convergence rate of $V_N^*$, but also characterizes the uncertainty about $\bbe_{P_N}[Z]$, which is the Bayesian average estimator of $\bbe_{P_*}[Z]$ for a general function $Z$. More specifically, the result in Proposition~\ref{pr-convergence} implies that for large $N$, the uncertainty in the Bayesian average estimator is inversely proportional to the (Fisher) information that the data $\xi^{(N)}$ carries about the true parameter $\theta^*$, and is proportional to the square of
$
\nabla_\theta \bbe_{P_{\theta^*}}[Z]=\nabla_\theta \int_\Xi
Z (\xi)  f(\xi|\theta^*)d\xi,
$
provided the integration and differentiation can be interchanged. 
%and where $P_\theta$ is the distribution defined by the pdf $f(\cdot|\theta)$.  
Note that $\nabla_\theta \bbe_{P_{\theta^*}}[Z]$ can be interpreted  as  the sensitivity of the function $\bbe_{|\theta}[Z]$  to the perturbation of the parameter value around $\theta^*$.

\subsection{Example: inventory control} \label{sec:inventory}
%\begin{example} [Inventory model]
In this subsection we apply the episodic Bayesian optimal control to the   classical inventory control problem and analytically show the the episodic value functions converge at the asymptotic convergence rate of $1/\sqrt{N}$.

{\rm
Consider the stationary inventory model
\[
\begin{array}{cll}
\min\limits_{u_t\ge 0} & \bbe_{P_*}\left[\sum_{t=1}^\infty\gamma^{t-1} \big(c u_t+\psi(x_t+u_t,D_t)\big)\right]\\
{\rm s.t.}& x_{t+1} = x_t+u_{t}-D_t,
\end{array}
\]
where $$\psi(y,d):=b[d-y]+h[y-d],$$
$c, b,h \in \bbr_+$  are the ordering cost, backorder penalty cost,
and holding cost per unit, respectively (with $b > c >0$), $x_t$ is the
current inventory level, $u_t$  is the order quantity, and $D_t\ge 0$  is
the demand at time $t$  which is a random iid process. The optimal policy is the base-stock  policy $\pi^*(x)=[x^*-x]_+$,  where $x^*=H^{-1}\left( \kappa\right)$ with
$$H(\nu):=P_*(D\le \nu)=\bbe_{P_*}[\ind_{(-\infty,\nu]}(D)]$$  being   the cdf of the demand,
 $\kappa:=\frac{b-(1-\gamma)c}{b+h}$
, and $\ind_A$ denotes the indicator function of set $A$  (e.g., \cite{zipkin}).  It follows that for $x\le x^*$ (cf., \cite{shacheng2021}),
 \begin{equation}\label{inval}
 V^*(x)=-c x+(1-\gamma)^{-1}\bbe_{P_*}\big[\gamma cD +(1-\gamma) c x^*+
\psi(x^*,D)\big].
 \end{equation}

Suppose that  the probability distribution of the random demand $D$  is modeled by pdf $f(\xi|\theta)$, $\theta\in \bbr$. Then $\pi_N(x)=[x_N-x]_+$,
where $x_N=H_N^{-1}\left( \kappa\right)$  with $H_N(\nu) =\bbe_{P_N}[\ind_{(-\infty,\nu]}(D)]$.
By  \eqref{zapprox} and \eqref{sigm}
we can use the following approximation
\[
H_N(\nu) \approx
H(\nu)+N^{-1/2}Y,
\]
where $Y\sim \N(0,\sigma^2)$,
with
\[
\sigma^2=\frac{1}{I(\theta^*)}\left(\int_{-\infty}^\nu \nabla_\theta f(\xi|\theta^*)  d\xi\right )^2=I(\theta^*)^{-1}\left( P_*\{S(\theta^*)\le \nu\}\right)^2.
\]
Suppose further  that the cdf  $H(\nu)$ has density $h(x^*)=  H'(x^*)>0$ at $\nu =x^*$. Then
\[
x_N  \approx
x^*-N^{-1/2} [h(x^*)]^{-1}Y.
\]
It follows that for $x\le \min\{x^*,x_N\}$ we have the following
\begin{equation}\label{piappr}
 \pi_N(x)-\pi^*(x)\approx N^{-1/2} [h(x^*)]^{-1}Y.
\end{equation}
Note that since the distribution of $Y$ is symmetrical around zero, the plus sign in the right hand side of \eqref{piappr} can be changed to the minus sign.

Also, by \eqref{inval} and \eqref{appr-2},
\begin{equation}\label{valap-1}
  V^*_N(x)-V^*(x)   \approx    (1-\gamma)^{-1}\bbe_{P_N-P_*}\big[\gamma cD +
\psi(x^*,D)\big].
\end{equation}
Together with \eqref{zapprox} this implies
\begin{equation}\label{valap-2}
  V^*_N(x)-V^*(x)   \approx    N^{-1/2} (1-\gamma)^{-1} W,
  \end{equation}
where $W$ has normal distribution with zero mean and variance
$$I(\theta^*)^{-1}(\bbe_{P_*}[(\gamma cD +
\psi(x^*,D)) S(\theta^*)])^2.$$

}%$\hfill\square$

%{\color{blue}
%\subsection{Regret analysis}
%Since in each episode the policy $\pi_N$ is applied to the true system and incurs a cost, we also consider the regret associated with this sequence of policies $\{\pi_N, N=1,2, \ldots \}$. We define the regret over $T$ episodes as follows:
%$$
%R(T) = \sum_{N=1}^T \bbe_{P_*} \left[c(x_N, \pi_N(x_N),\xi_N) - c(x_N, \pi^*(x_N), \xi_N)\right],
%$$
%where $\pi^*$ is the optimal policy of the true problem \eqref{exp-1}.
%Another definition of regret in the literature of learning discounted MDP is
%$$
%R(T) = \sum_{N=1}^T \bbe_{P_*} \left[V_N(x_N) - V^*(x_N)\right].
%$$
%
%First, we assume that the stage-wise cost function is Lipschitz in the control variable.
%\begin{assumption}
%There exists a constant $L_c$ such that for all $x\in\mathcal{X}$ and $\xi\in \Xi$,
%$$
%|c(x, u, \xi) - c(x, u', \xi)| \leq L_c\|u-u'\|, \forall u, u'\in \mathcal{U}.
%$$
%\end{assumption}
%
%For each $N$, we have
%\begin{align*}
%\bbe_{P_*} \left[c(x_N, \pi_N(x_N),\xi_N) - c(x_N, \pi^*(x_N), \xi_N)\right] &\leq  L_c\|\pi_N(x_N) - \pi^*(x_N)\|
%\end{align*}
%
%If we can show the convergence rate of $\pi_N$ to $\pi^*$, then we can show the regret grows sublinearly in $T$.}

%\setcounter{equation}{0}
\section{Computational Method}
\label{sec-comp}

As mentioned in Section~\ref{sec-boc}, the proposed episodic Bayesian optimal control requires solving a Bayesian average Bellman equation in every episode. In this section we look into a class of problems with a convex structure and develop an efficient computational method. In a nutshell, we first approximate the Bellman equation with a sample average and then extend the Stochastic Dual Dynamic Programming (SDDP) approach to the episodic Bayesian problem. SDDP is a dynamic cutting plane method that was originally developed for solving multi-stage linear stochastic programming problems  \cite{per1991}. We  extend SDDP to the infinite-horizon and episodic setting. It approximates the optimal value function by cutting planes, and thus it  produces a lower bound on the (convex) optimal value function. An upper bound can also be constructed by estimating value of the obtained (feasible) policy using methods such as  Monte Carlo simulation. Therefore, SDDP provides a (statistical) guarantee for accuracy of the obtained policy even in the case of continuous state and action spaces.  Please note that for the development of the computational method, the assumptions in Section~\ref{sec:convergence} are not required and new assumptions are introduced in this section.

\subsection{Sample average approximation of Bellman equation}
\label{sec-disc}

In order to solve Bellman equation \eqref{bellm-3} numerically the distribution $P_N$ should be discretized. We approach such discretization by generating a random sample using Monte Carlo sampling techniques.
There are two somewhat natural ways in which a random sample from $P_N$ can be generated. One straightforward way is to generate a random sample $\xi^1,...,\xi^M$,  of size $M$ from the distribution $P_N$. That is,     a point $\theta$ is generated from the posterior distribution   $p(\theta|\bxi^{(N)})$, and then conditional on $\theta$ a random point $\xi$ is generated from the density $f(\xi|\theta)$. This procedure is repeated $M$ times, independently from each other, to generate the random sample  $\xi^1,...,\xi^M$. Another approach is to generate a  sample  $\theta^1,...,\theta^{K_1}$  from  the
 the posterior distribution of
$p(\theta|\bxi^{(N)})$, and conditional on $\theta=\theta^i$, $i=1,...,K_1$,  to generate  a random sample $\{\xi^{i1},...,\xi^{iK_2}\}$   from the  pdf $f(\xi|\theta)$. The total sample size then is $K_1K_2$. The first method can be viewed as a special case of the second method with  $K_1:=M$ and $K_2:=1$.

To estimate $\bbe_{\theta_N} \left[\bbe_{\xi|\theta}[Z(\xi)]\right]$, the sample average estimator is
$$
\hat{Z}_{K_1,K_2} = \frac{1}{K_1K_2}\sum_{i=1}^{K_1}\sum_{j=1}^{K_2}Z(\xi^{ij}).
$$
This estimator is unbiased, i.e.,
$\bbe_{P_N}[\hat{Z}_{K_1,K_2}]
$ is equal to the expectation with respect to the joint distribution of $\xi$ and the posterior distribution of $\theta$ given by the right-hand side of \eqref{post-2a}.
Thus, it is sufficient to analyze its variance
$$
\var(\hat{Z}_{K_1,K_2}) = \frac{1}{(K_1K_2)^2}\sum_{i=1}^{K_1}\sum_{j=1}^{K_2}\var(Z(\xi^{ij})).
$$
Recall that variance of a random variable can be written as 
\[
\var(X)=\bbe[\var(X|Y)]+\var[\bbe(X|Y)],
\]
where
$
\var(X|Y)=\bbe\big[\big(X-\bbe(X|Y)\big)^2  |Y\big]
$
is the conditional variance of $X$ given $Y$ and
$
\var[\bbe(X|Y)]=\bbe\big [(\bbe(X|Y)-\bbe(X))^2\big].
$
 Thus
 \begin{equation}\label{variance}
    \var[Z(\xi^{ij})]=\bbe[ \var(Z(\xi^{ij})|\theta^i)]+\var[\bbe(Z(\xi^{ij})|\theta^i)].
 \end{equation}
Note that both terms in the right-hand side of \eqref{variance} are the same for all $i=1,\ldots, K_1$ and all $j=1,\ldots,K_2$. Therefore,
the variance of the estimator $\hat{Z}_{K_1,K_2}$  only depends on $K_1K_2$,  the total sample size, and is not affected by the individual choice of $K_1$ and $K_2$. For convenience, we will use the first method above to generate samples.

So suppose that random sample
$\xi^1,...,\xi^M$, is generated from $P_N$  by the first method (we use the same sample size $M$ for all $N$).
Then Bellman equation \eqref{bellm-3} is discretized to
\begin{equation}
 \label{disk-1}
\V_{N}(x)= \frac{1}{M}   \inf\limits_{u\in \U} \sum_{j=1}^M
   \big [c(x,u,\xi^j)+\gamma
\V_{N} \big(F(x,u,\xi^j)  \big)\big].
\end{equation}
Let
\[
\T_{N}(V)(x):=\frac{1}{M}   \inf\limits_{u\in \U} \sum_{j=1}^M
   \big [c(x,u,\xi^j)+\gamma
V \big(F(x,u,\xi^j)  \big)\big]
\]
be the corresponding  Bellman operator and $\V^*_{N}$ be  the  solution of equation \eqref{disk-1}. That is  $\V^*_{N}= \T_{N}(\V^*_{N})$.
Similar to \eqref{eq-est} we have that
\begin{equation}\label{disk-m}
 \|V^*_N-\V^*_{N}(x)\|_\infty\le (1-\gamma)^{-1}\|\cT_N(V^*_N)-\T_{N}(V^*_{N})\|_\infty.
\end{equation}
Under mild regularity conditions, for given    $\e>0$, $N$  and $\gamma\in (0,1)$ the
probability     $$P_N\big\{\|\cT_N(V^*_N)-\T_{N}(V^*_{N})\|_\infty\ge \e\big\}$$ converges to zero exponentially fast with the increase of the sample size $M$  (cf.,  \cite[section 8.4.2]{SDR}). By \eqref{disk-m} this implies
that   the
probability   $P_N\big\{\|V^*_N-\V^*_{N}\|_\infty\ge \e\big\}$ converges to zero exponentially fast with the increase of the sample size $M$.
% A Central Limit Theorem type results for $\V^*_{N}(x)$, considered as an estimate of $V^*_N(x)$,  are given in \cite{shacheng2021}{\color{blue}, whose main conclusion is that the sample size $M$ needed to maintain the relative stochastic error of the discretized problem is not sensitive to the discount factor $\gamma$ even if $\gamma$ is very close to one}.

%We  would like to show that  $\bar{\V}_N$   converges in some probabilistic sense to $V^*_N$ as $M\to\infty$. Consider operator $\cT_N$ defined in \eqref{operat-2}. Since $V^*_N$ is the fixed point of we have that it is the limit, in the norm topology of $\bbb(\X)$,  of   sequence $v_k\in \bbb(\X)$ defined iteratively $v_{k+1}=\cT_N(v_k)$, $k=1,...,$ starting from any $v_1\in \bbb(\X)$. Similarly $\bar{\V}_N$ is the limit of sequence $w_{k+1}=\T_N(w_k)$, $k=1,...,$, starting from the same point $v_1$.

\subsection{SDDP for infinite-horizon optimal control}

In this subsection, we develop a SDDP procedure for solving the discretized Bellman equation \eqref{disk-1} of one fixed episode.
For the sake of simplicity
we assume now  that the set $\X=\bbr^n$.
We consider   a setting where value functions, given by solutions of Bellman equations, are convex.
In order to enforce such convexity we make the following assumption.

\begin{assumption}
\label{ass-conv}
The  set  $\U$ is  convex, for any $\xi\in \Xi$  the
cost function $c(x,u,\xi)$ is convex in $(x,u)\in \bbr^n\times \U$,   and the mapping  $F(x,u,\xi)$ is affine, i.e.,
\begin{equation}\label{affin}
 F(x,u,\xi)=A(\xi) x+ B(\xi) u+ b(\xi).
\end{equation}
\end{assumption}
Under Assumption \ref{ass-conv}, value functions, given by solution of equations \eqref{bellm-1}, \eqref{bellm-3},  \eqref{disk-1},  are convex. For state-dependent controls, the convexity still holds (for details, please see \ref{EC:convexity} in the electronic companion.) %In that case cutting plane type algorithms can be applied to solve the discretized dynamic programming equation \eqref{disk-1}.

We discuss now a cutting planes algorithm for solving the infinite horizon {\em discretized} problem \eqref{disk-1}.
 We say that an affine function
$\ell(x)=\alpha^\top x+\beta$ is a {\em cutting} plane of a convex function  $\V:\bbr^n\to \bbr$ if $\V(x)\ge \ell(x)$ for all $x\in \bbr^n$. We say that a cutting plane $\ell(x)$ is a {\em supporting} plane of $\V(x)$,  at a point $\bar{x}\in \bbr^n$,  if $\V(\bar{x})=\ell(\bar{x})$.  A supporting plane of $\V(x)$  at   $\bar{x}$ is given by affine function $\ell(x)=\V(\bar{x})+g^\top (x-\bar{x})$, where $g$ is a subgradient of $\V(x)$ at $\bar{x}$.
Note that gradient of affine function
$\ell(x)=\alpha^\top x+\beta$ is $\nabla \ell(x)=\alpha$ for any $x$.

The cutting planes algorithm, of the SDDP type, approximates  value function $\V^*_N(x)$ by its cutting planes.
The traditional SDDP, for finite-horizon stochastic programming, consists of two steps in each iteration: a forward step that generates a path of trial points, and a backward step that adds cutting planes at these trial points to approximate the value functions at each stage. Since the value function is stationary in the infinite-horizon problem considered here, we can simplify these two steps: the forward step only simulates the next state (which is a new trial point), and the backward step adds one cutting plane on this new trial point.
%  In the backward step of the algorithm additional cutting planes are added to the   approximation of the value function.

More specifically,  let  $\underline{\V}_N(x)$
 be the current approximation of $\V^*_N(x)$,  given by the maximum of a finite number of cutting planes. Note that by the definition of cutting planes,  we have that
$\V^*_N(\cdot)\ge \underline{\V}_N(\cdot)$. Given the current state (trial point) $\bar{x}\in \bbr^n$, compute the control
\begin{equation}
 \label{sddpgr-2}
 \bar{u}\in \arg\min_{u\in \U}   \sum_{j=1}^M  \big [c_j (\bar{x},u)+\gamma
\underline{\V}_N \big(A_{j}\bar{x} +B_{j}u +b_{j}  \big)\big],
\end{equation}
where $c_j(x,u):=c(x,u,\xi^j)$, $A_j:=A(\xi^j)$, $B_j:=B(\xi^j)$, $b_j:=b(\xi^j)$.
Then, additional  cutting plane is computed as $\ell(x)= \underline{v}_N+g^\top (x-\bar{x})$, where
\begin{eqnarray}
 \underline{v}_N &=&  M^{-1}   \sum_{j=1}^M  \big [c_{j}(\bar{x},\bar{u})+\gamma
\underline{\V}_N \big(A_{j}\bar{x} +B_{j}\bar{u} +b_{j}  \big)\big],\\
g &=& M^{-1}   \sum_{j=1}^M
   \big [\nabla_x c_{j}(\bar{x},\bar{u})+\gamma A_{j}^\top \nabla
\underline{\V}_N \big(A_{j}\bar{x} +B_{j}\bar{u} +b_{j}  \big)\big].
\label{sddpgr-1}
\end{eqnarray}

 The   subgradient   $\nabla \underline{\V}_N(x')$ at the point $x':=A_{j}\bar{x} +B_{j}\bar{u} +b_{j}$ is computed as follows.
Let  $\underline{\V}_N(x)=\max_{\eta\in \I}\ell_\eta(x)$ be the  current representation of  $\underline{\V}_N(x)$
by cutting planes (affine functions)  $\ell_\eta(x)$, $\eta\in \I$.
 Let $\hat{\eta}\in \I$ be such that $\underline{\V}_N(x')=\ell_{\hat{\eta}}(x')$, i.e., $\ell_{\hat{\eta}}(x)$ is the supporting plane of $\underline{\V}_N(x)$ at $x'$.
Then $\nabla \underline{\V}_N(x')=\nabla\ell_{\hat{\eta}}(x')$.
The forward step generates the trial point by simulating the next state from the current state and control: $x':=A_{j}\bar{x} +B_{j}\bar{u} +b_{j}$, where $A_{j}$, $B_j$, and $b_j$ are computed on a randomly sampled $\xi_j$. For more implementation details on how to compute the subgradient, please refer to \ref{EC:subgradient} in the electronic companion.

%The forward step of the algorithm has two functions. It constructs a (statistical) upper bound for the optimal value of the considered    discretised  problem,   and generates trial points for the subsequent iteration of the backward step.
%  The constructed approximation $\underline{\V}_N(x)$  of the value function defines a feasible policy for the corresponding infinite horizon problem. Therefore value of this policy gives an upper bound for the opimal value of the discretized problem.

The approximation  $\underline{\V}_N(x)$ of value functions defines a policy, and the value of the policy can be estimated as follows.
 Let $\hat{\xi}_t$, $t=1,...$, be a sample path (an iid sequence)  from the discretized distribution. Such  sample path is generated by taking at every stage  one of the values   $\xi^{j}$,  $j=1,...,M$, at random    with equal  probability $M^{-1}$.
 Starting with the initial value $\hat{x}_1=x_1$, the  state and control variables  are computed sequentially going forward in time. At stage $t=1,...,$  given   value  $\hat{x}_t$ of the state vector,  the corresponding control is computed as
 \begin{equation}
 \label{sddpgr-2a}
 \hat{u}_t\in \arg\min_{u\in \U}   \sum_{j=1}^M  \big [c_{j}(\hat{x}_t,u)+\gamma
\underline{\V}_N \big(A_{j}\hat{x}_t +B_{j}u +b_{j}  \big)\big].
\end{equation}
 Then next value $\hat{x}_{t+1}=A(\hat{\xi}_t)\hat{x}_t+B(\hat{\xi}_t)\hat{u}_t+
b(\hat{\xi}_t)$  of the state vector is computed. Consequently  with this sample path $\{\hat{x}_t,\hat{u}_t\}$,   an unbiased point estimate of  the   value  of the constructed  policy is
\begin{equation}\label{infinit-series}
\sum_{t=1}^\infty   \gamma^{t-1} c (\hat{x}_t, \hat{u}_t,\hat{\xi}_t).
\end{equation}
The expected value of \eqref{infinit-series} gives an upper bound for the optimal value  $\V^*_{N}(x_1)$ of the discretized Bellman equation \eqref{disk-1}. The tightness of this upper bound depends on the accuracy  of the constructed policy.

%This evaluation of the policy also provides a (statistical) upper bound on the optimal value function $\V^*_{N}$ of the discretized Bellman equation \eqref{disk-1}. This upper bound is only statistical instead of being exact, due to the simulation error in the evaluation.

In order to compute the infinite sum \eqref{infinit-series},  it should be truncated,  i.e.,  it is approximated by a finite sum $\sum_{t=1}^T  \gamma^{t-1} c (\hat{x}_t, \hat{u}_t,\hat{\xi}_t)$.
To determine the truncated length, let $\kappa$ be  such that   $|c(x,u,\xi)|\le \kappa$ for all considered $(x,u,\xi)$.
 Then  the error of the  truncation can be bounded  as
 \[
 \begin{array}{ll}
 \left|\sum_{t=1}^\infty \gamma^{t-1}c (\hat{x}_t, \hat{u}_t,\hat{\xi}_t)
 -\sum_{t=1}^T  \gamma^{t-1}c (\hat{x}_t, \hat{u}_t,\hat{\xi}_t)\right|
 \le (1-\gamma)^{-1}\gamma^T \kappa.
 \end{array}
 \]
 Consequently, the time horizon $T$ is chosen in such a way that the truncation error is less than the prescribed tolerance $\epsilon$, i.e.,
\begin{equation}\label{truncation}
T\geq \log_{\gamma}\frac{\epsilon(1-\gamma)}{\kappa}.
\end{equation}

\subsection{Extension of SDDP to episodic setting}
The SDDP algorithm described above is for one episode, i.e., for a fixed $N$, where the posterior pdf  $ p(\theta|\bxi^{(N)})$ is fixed. It could be extended to the multi-episode setting in a naive way by restarting from the scratch for every $N$. However, we can reuse value function approximation from previous episode to warm start the SDDP procedure for solving the problem at the next episode.
To initialize the SDDP procedure, a lower bound on the value function is needed. Hence, a natural question is whether a lower approximation $\underline{\V}_N$ obtained from episode $N$ is still a lower bound for  $\V_{N+1}^*$, the optimal value function of discretized problem at episode $N+1$. 

Let us note that if $V^*$ is the solution of Bellman equation \eqref{bellm-1} and $\cT$ is the corresponding Bellman operator, then $V \geq \cT (V)$ implies that   $V \geq V^{*}$. Indeed it follows then by monotonicity of $\cT$ that $V \geq \cT (V)\ge \cT^2(V)\ge \cdots \ge \cT^k(V)$, and hence the assertion follows  since $\|\cT^k(V)-V^*\|_\infty$  tends to zero as $k\to \infty$. And similarly if $V \leq \cT (V)$, then $V \leq V^{*}$. That is, the sub-solution and super-solution of the optimality equation serve as lower and upper bound for the optimal value function.
Therefore 
for $\underline{\V}_{N} \leq \V_{N+1}^*$ to hold, it suffices to check 
$$
\underline{\V}_{N}(x) \leq \T_{N+1} (\underline{\V}_{N})(x),
$$
where $\T_{N+1}$
is the Bellman operator of the discretized problem at episode $N+1$, which can be rewritten as
\begin{align}\label{Bellman-LR}
\T_{N+1}(\V) (x)= \inf_{u\in \mathcal{U}}\frac{1}{M}\sum_{j=1}^M
\left[\frac{p(\theta^j|\bxi^{(N+1)})}{p(\theta^j|\bxi^{(N)})}\left[c(x,u,\xi^j)+\gamma \V(F(x,u,\xi^j))\right] \right].
\end{align}
Since the likelihood ratio $p(\theta|\bxi^{(N+1)})/p(\theta|\bxi^{(N)})$ could be smaller than 1,  it is easy to see that $\underline{\V}_{N}(x) \leq \T_{N+1} (\underline{\V}_{N})(x)$ does not necessarily hold. However, since $p(\theta|\bxi^{(N)})$ converges as $N$ goes to infinity, we can expect that the likelihood gets close to 1 when $N$ is large and hence the lower approximation $\underline{\V}_{N}(x)$ would also be a close lower approximation for the next episode.

Motivated by the above observation, the next question is whether we can reuse some of the cutting planes in $\underline{\V}_N$ to form a lower bound on the value function of next episode.  To answer this question, we take two steps: first, we identify the cutting planes that are individually a lower bound, and we then show that the maximum of these cutting planes forms  a valid lower bound for the value function of next episode. As mentioned above, to identify whether a cutting plane in $\underline{\V}_{N}$, denoted by $\ell_{N, \eta}(x)$, is a lower bound of $\V_{N+1}^*$, it suffices to check the following condition:
\begin{eqnarray}\label{cut}
\ell_{N,\eta}(x) &\leq& \T_{N+1} (\ell_{N,\eta})(x) \nonumber\\
&=& \inf_{u\in \mathcal{U}}\frac{1}{M}\sum_{j=1}^M
\left[\frac{p(\theta^j|\bxi^{(N+1)})}{p(\theta^j|\bxi^{(N)})}\big[c(x,u,\xi^j)
+\gamma \ell_{N,\eta}(F(x,u,\xi^j)\big] \right].
\end{eqnarray}
%For computation, we replace the exact Bellman operator $\cT_{N+1}$ by the discretized approximation:
%\begin{equation} \label{cut}
%\ell_{N,\eta}(x) \leq ,
%\end{equation}
%where $(\theta_i,\xi_i)$ are the samples  at episode $N$, i.e., drawn from $p_N(\theta)p(\xi|\theta)$.
That is equivalent to checking the following holds:
\begin{equation}\label{cut_check}
\min_{x\in\mathcal{X},u\in \mathcal{U}} \frac{1}{M}\sum_{j=1}^M
\left[\frac{p(\theta^j|\bxi^{(N+1)})}{p(\theta^j|\bxi^{(N)})}\big[c(x,u,\xi^j)
+\gamma \ell_{N,\eta}(F(x,u,\xi^j)\big] \right] 
- \ell_{N,\eta}(x) \geq 0.  
\end{equation}

If the sets $\mathcal{X}$ and $\mathcal{U}$ are polyhedral, the function $c(x,u,\xi)$ is piecewise linear convex, and F is affine as in \eqref{affin}, then \eqref{cut_check} is a linear program (LP). This condition is checked once for every cutting plane in the lower approximation of the previous episode, and hence the overhead incurred in each episode is the time of solving the LP \eqref{cut_check} multiplied by the total number of cutting planes. Note that the SDDP procedure carries out multiple iterations in each episode, and each iteration requires solving \eqref{sddpgr-2a}, which is equivalent to solving an LP whose number of constraints is equal to the total number of cutting planes in the current lower approximation. Therefore, the computational time of checking condition \eqref{cut_check} is only a small fraction of the time used by the SDDP procedure. Moreover, warm-start by reusing previous cutting planes usually reduces the number of iterations for SDDP to converge, and hence the overhead incurred by checking condition \eqref{cut_check} is usually offset by the saved computational time of the SDDP procedure. This is also empirically observed in our numerical experiments in Section~\ref{sec:numerical}.
 
Now suppose we have found multiple cutting planes that are valid (i.e., the inequality \eqref{cut} holds), we can show that they are also jointly valid, that is, the lower approximation constructed from maximum of these cutting planes also satisfies the inequality. Suppose there are two valid cutting planes, denoted as $a_1 x + b_1$ and $a_2 x + b_2$ respectively. Then 
\begin{align*}
    & \min_{u \in \mathcal{U}} \sum_{j=1}^{M} \frac{p(\theta^j|\bxi^{(N+1)})}{p(\theta^j|\bxi^{(N)})}[c(x,u,\xi^j)+\gamma \max \{a_1 F(x,u,\xi^j) + b_1, a_2 F(x,u,\xi^j) + b_2\}] \\
    = & \min_{u \in \mathcal{U}} \max_{i \in \{1,2\}} f_i(u)  \;\; \Big(\text{where} \; f_i(u)=\sum_{j=1}^{M}\frac{p(\theta^j|\bxi^{(N+1)})}{p(\theta^j|\bxi^{(N)})}[c(x,u,\xi^j)+\gamma(a_i F(x,u,\xi^j) + b_i)]\Big)\\
    = & \min_{u \in \mathcal{U}} \max_{i \in \{1,2\}} f(i,u,u)  \quad (\text{where} \; f(i,u_1,u_2)=f_1(u_1)  \text{ if } i = 1, \text{and } f_2(u_2) \text{ if }  i = 2)\\
    \geq & \max_{i \in \{1,2\}} \min_{u \in \mathcal{U}} f(i, u, u) 
    \geq   \max_{i \in \{1,2\}} \min_{u_1 \in \mathcal{U}, u_2 \in \mathcal{U}} f(i, u_1, u_2)\\
    = & \max \{\min_{u_1 \in \mathcal{U}}f_1(u_1), \min_{u_2 \in \mathcal{U}}f_2(u_2)\} 
    \geq   \max \{a_1 x + b_1, a_2 x + b_2\}.
\end{align*}
By induction, we can aggregate all the cutting planes from $\underline{\V}_{N}$ that are lower bounds for episode $N+1$ and use that approximation as a warm start.

Putting the ideas above together, we present the algorithm for episodic Bayesian optimal control with SDDP in Algorithm~\ref{algorithm: episodic_SDDP}. In a nutshell, in each episode we adopt the SDDP procedure to solve the Bayesian average estimate of the stochastic optimal control problem with a fixed posterior distribution, and then exercise the obtained policy to transition to the next state;  with the observed additional data, we update the posterior distribution and warm start SDDP for the next episode.

% \begin{minipage}{1\linewidth}
\begin{algorithm}[th]
\SetAlgoLined
\SetKwInOut{Input}{input}\SetKwInOut{Output}{output}
\Input{initial state $x_1$; posterior distribution $p(\theta|\bxi^{(1)})$; SDDP iteration number $K$; sample size $M$.} %upper bound for cost function $\kappa$; discount factor $\gamma$; tolerance for the the error of the truncation $\cT$; initial lower approximation $\underline{V}_1(\cdot)$;
\Output{approximate optimal value functions $\underline{\V}_N, N=1, 2, \ldots$}
\For{$N =1, 2,\ldots$}{
\textbf{SDDP:}
%Observe $D_n$ data points of $\xi$ and update the posterior distribution $p_n$.
Set $\bar{x}_0 = x_N$.
If $N=1$, set initial lower approximation $\underline{\V}_1^1$; if $N\neq 1$, construct  $\underline{\V}_N^{1}$ using cutting planes from $\underline{\V}_{N-1}$ that satisfy \eqref{cut_check}.\\

\For{$k \leftarrow 1$ \KwTo $K$}{
Draw iid samples $\xi^1,\ldots,\xi^M$, and compute the control, value, and subgradient:
\begin{align*}
    \bar{u} \in \frac{1}{M} &\arg\min_{u \in \mathcal{U}} \sum_{j=1}^{M} c_j(\bar{x}_{k-1},u,\xi^j) + \gamma \underline{\V}^{k}_n(A_j \bar{x}_{k-1}+B_j u + b_j),\\
  \underline{v}_k &= \frac{1}{M} \sum_{j=1}^{M} c_j(\bar{x}_{k-1},\bar{u},\xi^j) + \gamma \underline{\V}^{k}_n(A_j \bar{x}_{k-1}+B_j \bar{u} + b_j), \\
   g &= \frac{1}{M} \sum_{j=1}^{M} \nabla_x c_j(\bar{x}_{k-1},\bar{u}, \xi^j) + \gamma A_j^{T} \nabla \underline{\mathcal{V}}^{k}_N(A_j \bar{x}_{k-1} +B_j \bar{u} + b_j).
\end{align*}\\
Update the lower approximation: $\underline{\V}_N^{k+1} = \max\{\underline{v}_k+g^{T}(x-\bar{x}_k), \underline{\V}_N^{k}\}$.\\
Generate the next trial point: randomly draw a sample $\xi^j$ and set $\bar{x}_{k} = A_j\bar{x}_{k-1}+B_j\bar{u}+b_j.$
}
%\textbf{Lower approximation:} Set $\underline{V}_{n}=\underline{\mathcal{V}}_{n}^{K}$. \\
\textbf{Policy Evaluation (optional):} Compute the truncated time horizon $T$ by \eqref{truncation}. Draw a sample path $\{\hat{\xi}_1,\ldots,\hat{\xi}_T\}$, and compute $$\hat{\V}_N = \sum_{t=1}^T \gamma^t c(\hat{x}_t, \pi_N(\hat{x}_t), \hat{\xi}_t),$$ where $\pi_N$ is the policy defined by $\underline{\V}_{N}$. Note that $\hat{\V}_N$ provides a statistical upper bound on the episodic optimal value function (of the truncated problem).}
\textbf{Updating:} Exercise the policy $\pi_N$ to transition to the next state $x_{N+1} = F(x_N, \pi_N(x_N), \xi_{N+1})$; update the posterior to $p(\theta|\bxi^{(N+1)})$ with $\xi_{N+1}$.

\caption{Episodic Bayesian optimal control with SDDP}
\label{algorithm: episodic_SDDP}
\end{algorithm}
% \end{minipage}

\section{Numerical Results}\label{sec:numerical}

In this section, we carry out numerical experiments to verify the theoretical results and demonstrate the performance of Algorithm~\ref{algorithm: episodic_SDDP}  on the inventory control problem that is stated in Section~\ref{sec:inventory}.

\subsection{One-dimensional case} \label{sec: numerical one-dim}
We first test on the single-product (one-dimensional) inventory control problem. The parameter settings are as follows:
%, where %we aim to solve
%$$
%\begin{array}{ll}
%\min _{u_t \geq 0} & \mathbb{E}_{P_*}\left[\sum_{t=1}^{\infty} %\gamma^{t-1}\left(c u_t+\psi\left(x_t+u_t, D_t\right)\right)\right] \\
%\text { s.t. } & x_{t+1}=x_t+u_t-D_t,
%\end{array}
%$$
%where
%$$
%\psi(y, d):=b[d-y]+h[y-d]
%$$
 unit ordering cost $c = 1$, unit holding cost $h=2$, unit backorder penalty cost $ b=3 $, and customer demand $D_t $ at time $t$ follows an exponential distribution with unknown mean $\theta = 10$. Recall that the optimal policy is the basestock policy $\pi^*(x)=\left[x^*-x\right]_{+}$, where $x^*=H^{-1}(\kappa)$ with $H(\cdot)$ being the cdf of the demand and $\kappa:=\frac{b-(1-\gamma) c}{b+h}$. It follows that for $x \leq x^*$,
\begin{equation} \label{inventory_value}
V^*(x)=-c x+(1-\gamma)^{-1} \mathbb{E}_{P_*}\left[\gamma c D+(1-\gamma) c x^*+\psi\left(x^*, D\right)\right].
\end{equation}
With this analytical form, we compute the true optimal value function as a benchmark to compare with the value functions of the episodic Bayesian problems.

To make the Bayesian updating in closed form, we choose a Gamma prior $\Gamma(\alpha_0,\beta_0)$ with $\alpha_0=\beta_0=1$ on the unknown parameter $\frac{1}{\theta}$. Given $n$ observed data $\xi_1,\ldots,\xi_n$, the posterior is also a Gamma distribution with parameter $\alpha' = \alpha_0 + n$, $\beta' = \beta_0 + \sum_{i=1}^n \xi_i$. In episode $N$, we update the Bayesian posterior with the observed data and then solve the above inventory control problem with probability distribution $P_*$ in \eqref{inventory_value} replaced by posterior probability distribution $P_N$ for the episodic value function $V_N^*$. { We use a large number ($10^5$) of samples from $P_*$ or $P_N$ to approximate the expectation in \eqref{inventory_value} to obtain the true optimal value function $V^*$ or the episodic optimal value function $V_N^*$.}  We then plot the integrated gap between $V_N^*$ and $V^*$ with the gap computed as $\int_{\mathcal{X}} |V_N^*(x) - V^*(x)|d\mu^*(x)$, where $\mu^*$ is the stationary distribution (or occupancy measure of states) under true optimal policy $\pi^*$ of the original problem.
\begin{figure}[tb]
        \centering
        \includegraphics[width=0.48\textwidth]{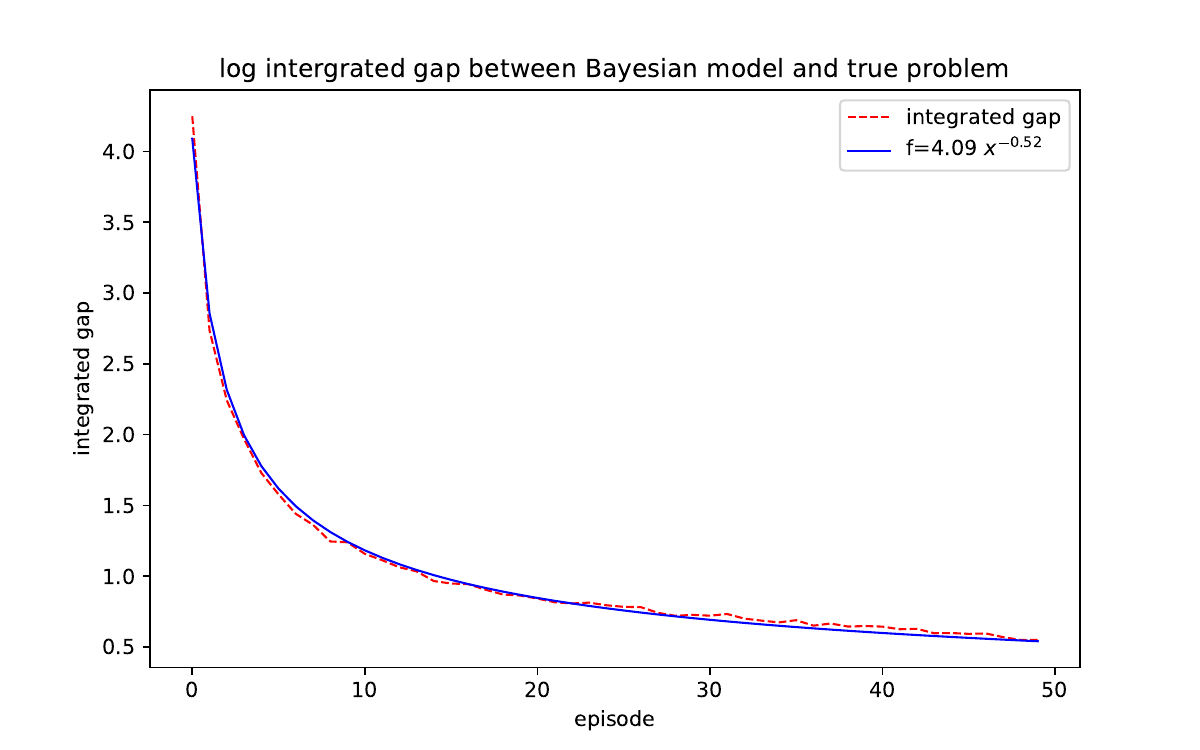}
       % \begin{subfigure}{0.48\textwidth}
       %     \includegraphics[width=1\textwidth]{numerical/normality_fit_discount factor gamma_0.6.pdf}
       %     \caption{Asymptotic normality}
       %     \label{fig: normality}
       % \end{subfigure}
\caption{Convergence rate of episodic value functions}
\label{fig: convergence rate}
\end{figure}

Figure \ref{fig: convergence rate} shows integrated gap between $V_N^*$ and $V^*$ in terms of episodes. The result is computed by running 200 replications. In each episode, new data of batch size $20$ is used to update the Bayesian posterior. To see how fast it converges to 0,  we fit a regression model  $f(x) = a * x^b$  with optimal values $a = 9.84$ and $b=-0.49$, which is shown as the blue curve in the figure. The close value of $b$ to $-0.5$ verifies the theoretical result of convergence rate $O(N^{-1/2})$ with $N$ being the number of episodes.

\begin{figure}[tb]
       \begin{subfigure}{0.48\textwidth}
           \includegraphics[width=1\textwidth]{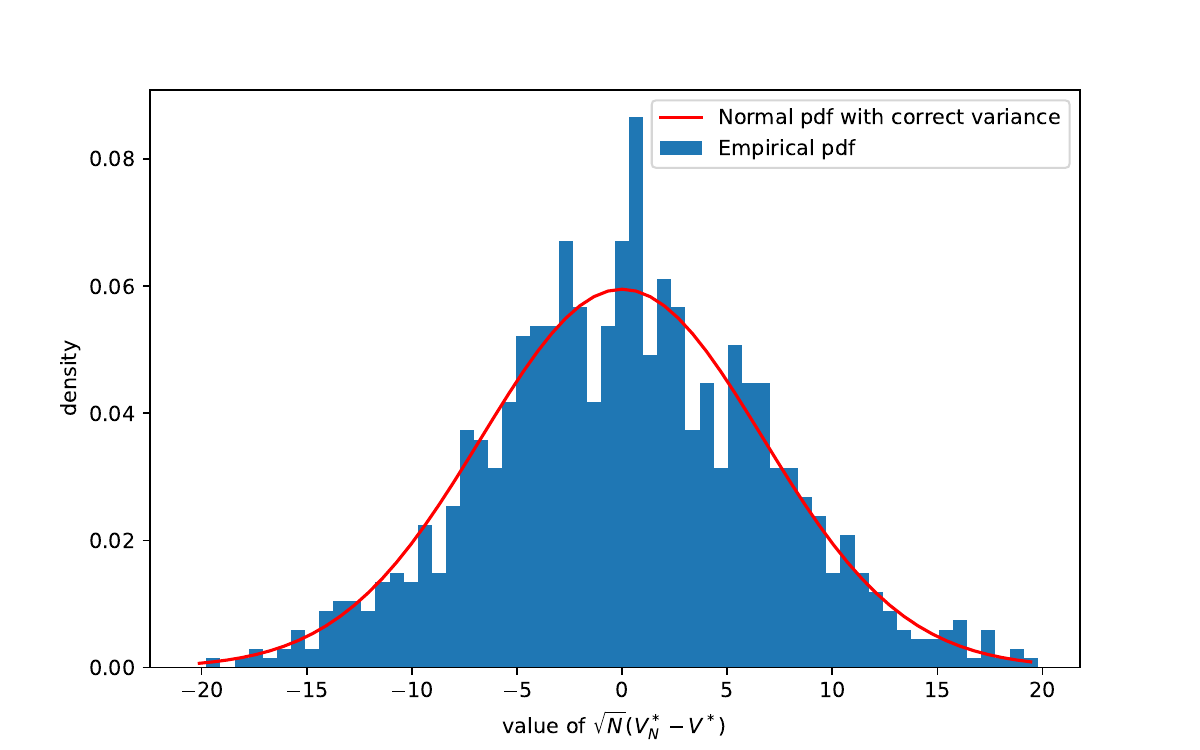}
           \caption{Empirical pdf}
           \label{fig: normality 1}
       \end{subfigure}
       \begin{subfigure}{0.48\textwidth}
           \includegraphics[width=0.85\textwidth]{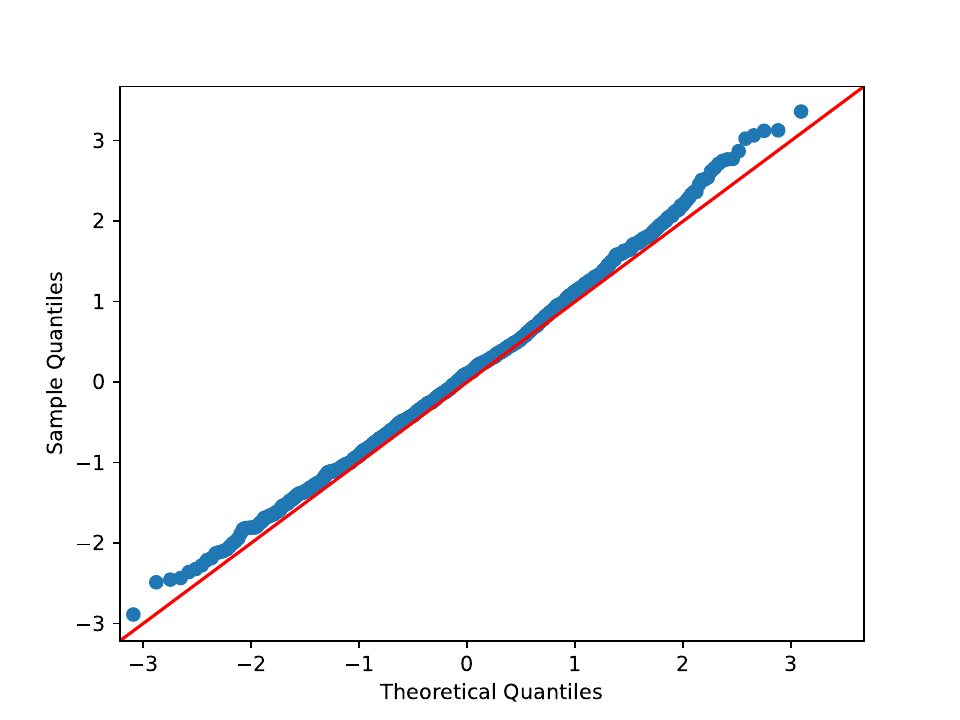}
           \caption{Normal QQ plot}
           \label{fig: normality 2}
       \end{subfigure}
\caption{Asymptotic normality of episodic value functions}
\label{fig: normality}
\end{figure}

In Figure \ref{fig: normality}, we test the asymptotic normality of $\sqrt{N}(V_N^* - V^*)$. For this specific test, we set the mean of demand $\theta = 1$ to reduce the Monte Carlo estimation error of $V_N^*$ in the experiment. Other parameters are the same as in Figure \ref{fig: convergence rate}. In Figure \ref{fig: normality 1}, the blue bar chart represents the empirical density of $\sqrt{N}(V_N^* - V^*)$ at $x_0 = 0$ and $N=100$, with $1000$ replications. The red curve represents the theoretical normal pdf with mean $0$ and standard deviation $\sigma^c$ computed according to \eqref{valap-2}. Figure \ref{fig: normality 2} is the QQ plot, which plots the quantile of the empirical distribution  $\frac{\sqrt{N}(V_N^* - V^*)}{\sigma^c}$ (empirical quantile) against the quantile of the standard normal distribution (theoretical quantile).  As both figures indicate, the empirical pdf coincides with the theoretical pdf well, verifying the asymptotic normality result in Section~\ref{sec-asymptotic}.

\begin{figure}[tb]
    {
        \centering
    \includegraphics[width=0.48\textwidth]{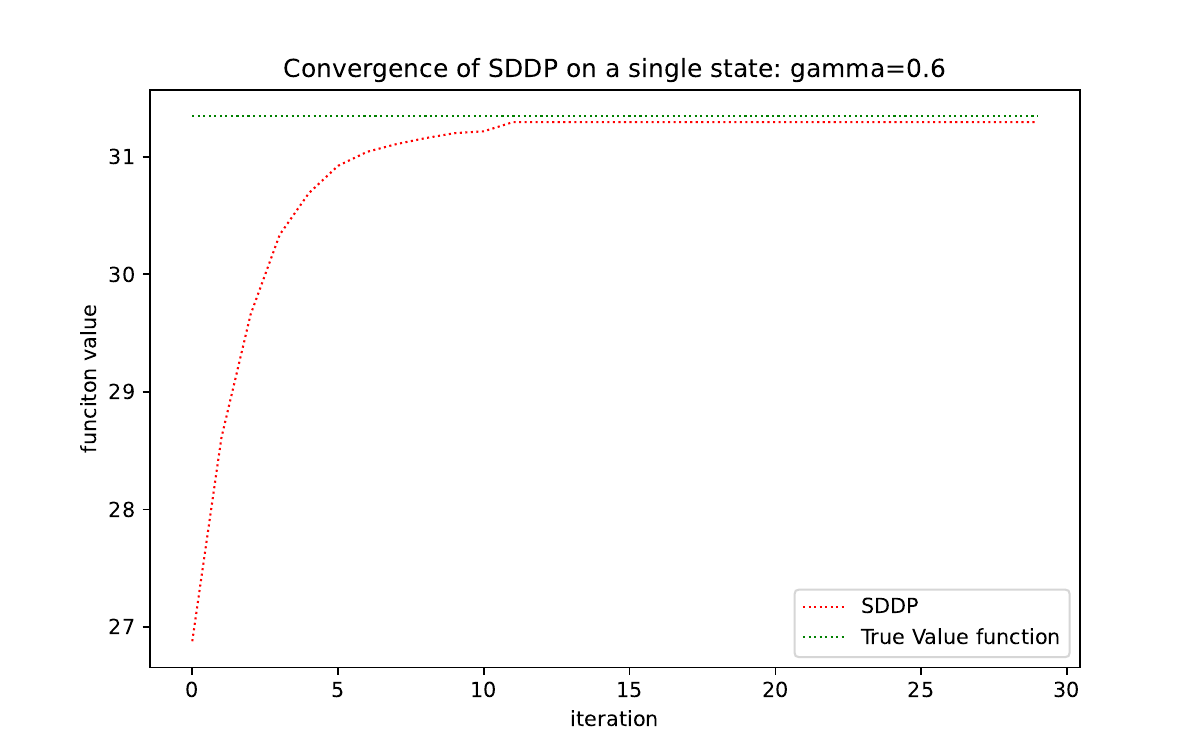}
        \includegraphics[width=0.48\textwidth]{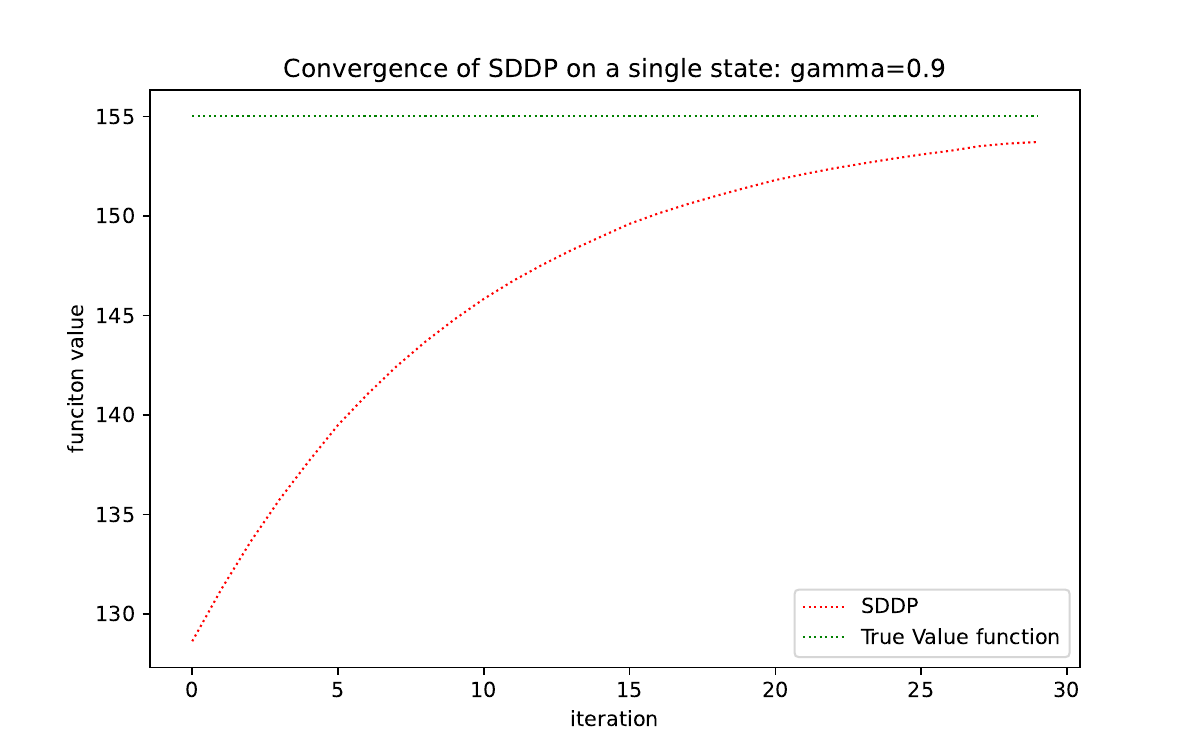}
    }
\caption{SDDP algorithm with single episode for 1-D inventory control.}
\label{fig: lower}
\end{figure}

Next, we demonstrate the performance of Algorithm~\ref{algorithm: episodic_SDDP} on this problem. We run two different versions of this algorithm, both of which use SDDP to solve episode-wise Bayesian problem but differs in whether warm starting with the lower approximation from the previous episode. { Both versions use the sample size $M = 100$.} In Figure \ref{fig: lower}, we plot the convergence of the SDDP algorithm on the final episode out of 5 episodes in total. The x-axis is the number of SDDP iteration number, and the y-axis is the function value (episodic value function $V_N^*$, and the approximate value function $\underline{\V}_N$ generated by Algorithm~\ref{algorithm: episodic_SDDP}) evaluated at the initial point $x_{0} = 1$. We vary the discount factor from 0.6 (left plot) to 0.9 (right plot).
Figure \ref{fig: lower} indicates the SDDP converges to the true value by maintaining an increasing lower approximation. It takes about 5 iterations to converge for $\gamma = 0.6$  and takes about 30 iterations to converge for $\gamma = 0.9$. The convergence is fast, considering the control and state space are both continuous.

\begin{figure}[tb]
    {
        \centering
        \includegraphics[width=0.48\textwidth]{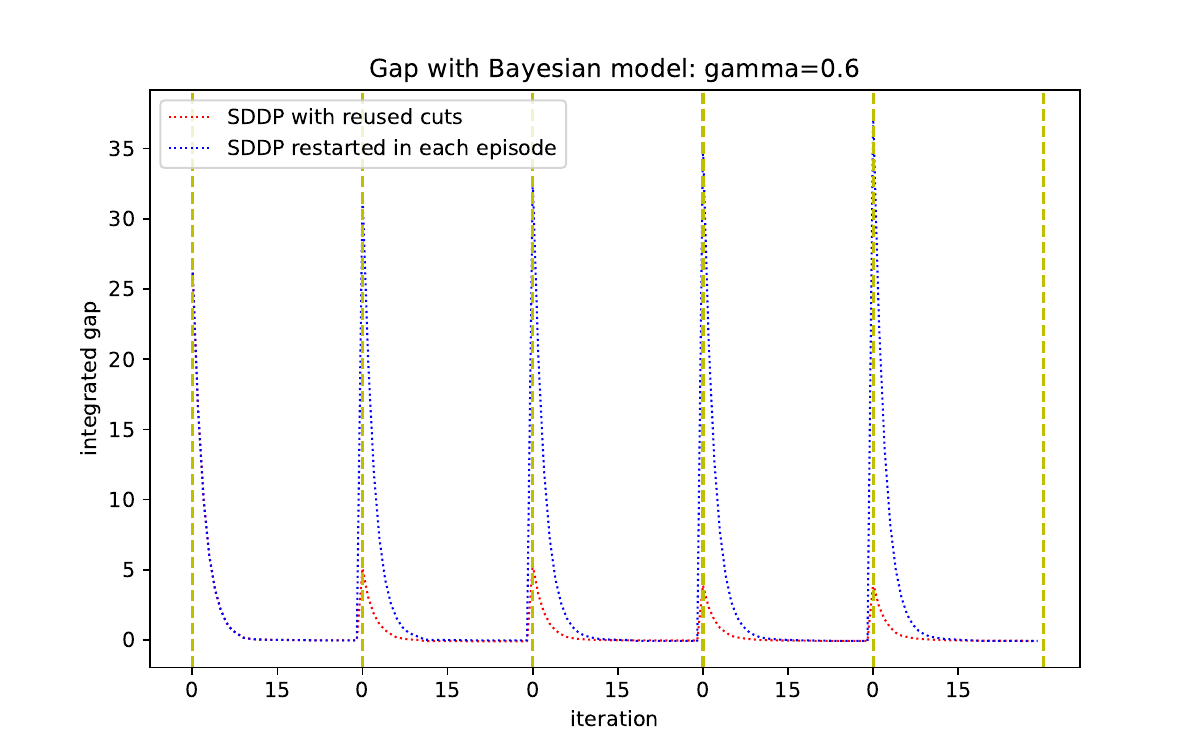}
        \includegraphics[width=0.48\textwidth]{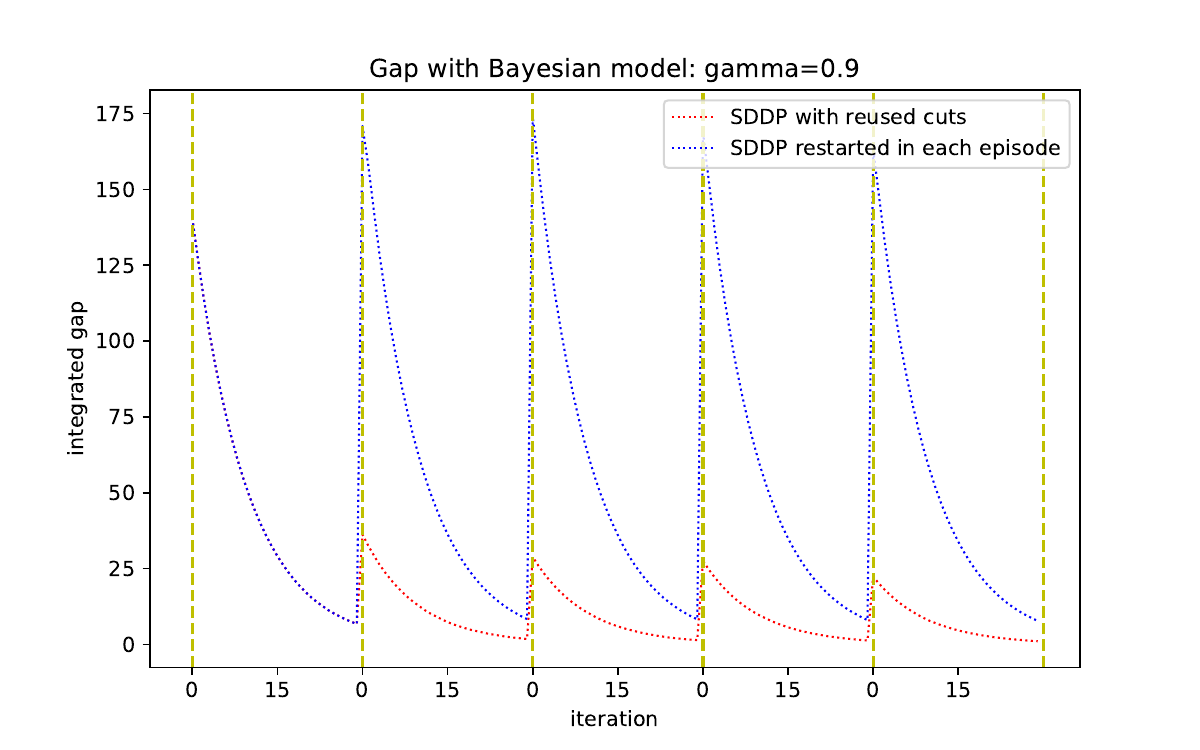}
    }
\caption{SDDP algorithm for 1-D Inventory Control.}
\label{fig: 2}
\end{figure}

 Figure \ref{fig: 2} plots the integrated gap between $\underline{V_N^*}$ and $V_N^*$ over
 %all potential states in terms of SDDP iterations, %using $L_1$ norm 
 %integrated with 
 the stationary distribution $\mu_N$ under the episodic optimal policy $\pi_N$, i.e., 
$ \int_\mathcal{X}(V_N^*(x) - \underline{\V}_N(x))d\mu_N(x)$, where we drop the absolute value since $\underline{\V}_N$ is always a lower bound of $V_N^*$.  
The vertical dashed line shows the switch between two successive episodes. The left and right plots vary in the discount factor $\gamma$.
Figure~\ref{fig: 2}  shows that within each episode SDDP converges to the optimal value $V_N^*$ efficiently. By comparing the algorithm with warm start (red curve) and without warm start (blue curve), we can see the benefits of reusing the previous lower approximation in the next episode. Without warm start, the lower approximation of each episode starts from a large value (about $35$ for $\gamma = 0.6$ and $150$ for $\gamma = 0.9$) since we initialize $\underline{\V}_N$ by using a loose lower bound. When reusing the previous lower approximation, the initial gap starts from a relatively small value (lower than $5$ for $\gamma = 0.6$ and $30$ for $\gamma = 0.9$).
%{\color{blue} Moreover, for $\gamma = 0.6$, with reused cuts, the gap between $\underline{V}_N$ and $V^*_N$ decreases to around $0$  within 5 iterations in the last 2 episodes, whereas without reused cuts, the gap converges to the same level after about 8 iterations. For $\gamma = 0.9$, the gap with reused cuts decreases to  around 1 at the end of the last episode while the gap without reused cuts only reaches around $8$.}
This is even more beneficial in the later episode as the Bayesian posterior concentrate more on the true value, and as a result, more than $50\%$ of cuts of the previous episode remain valid for the new episode.

\begin{figure}[tb]
        \centering
 \includegraphics[width=0.48\textwidth]{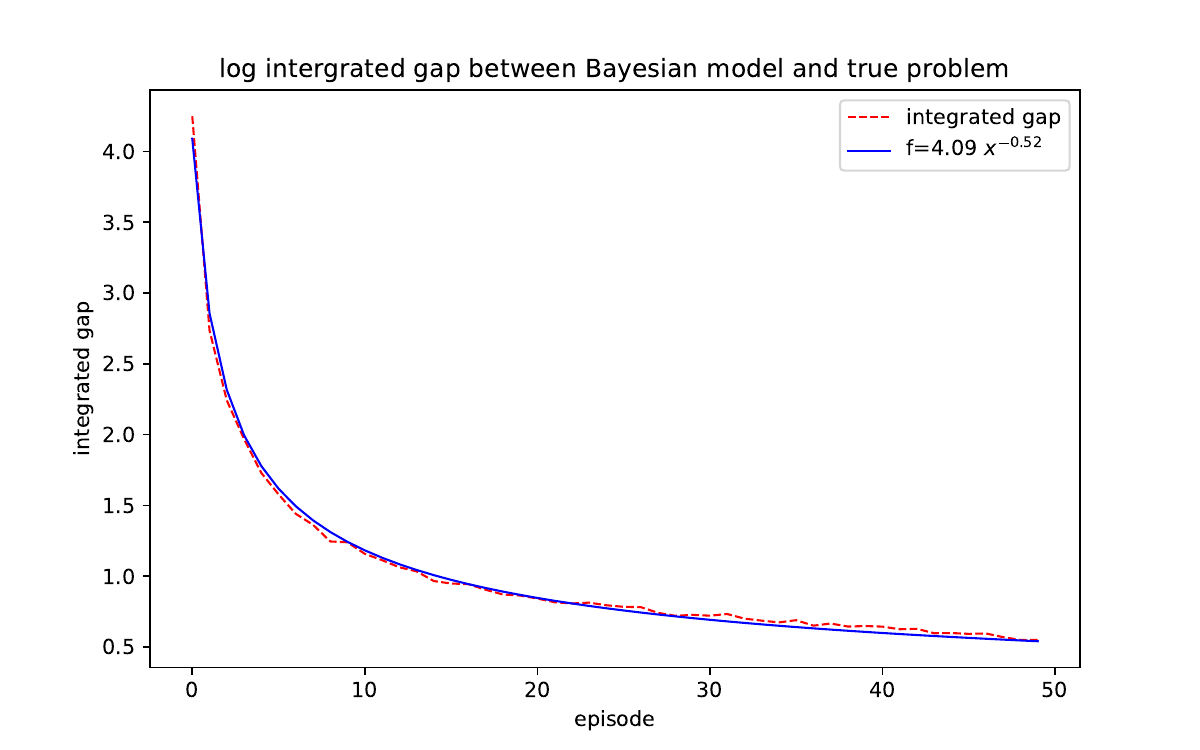}
\caption{Convergence rate of episodic value functions}
\label{fig: convergence multi}
\end{figure}

\begin{figure}[tb]
    {
        \centering
       \begin{subfigure}{0.48\textwidth}
           \includegraphics[width=1\textwidth]{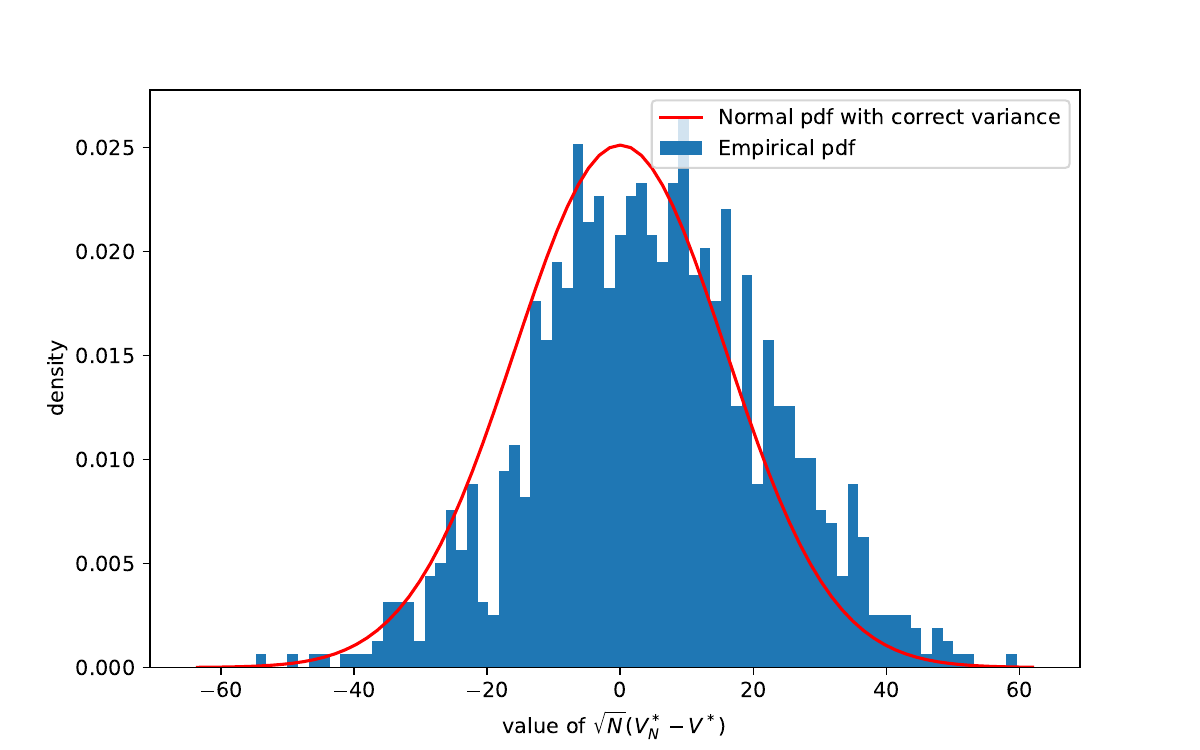}
           \caption{Empirical pdf}
           \label{fig: normality multi 1}
       \end{subfigure}
       \begin{subfigure}{0.48\textwidth}
           \includegraphics[width=0.85\textwidth]{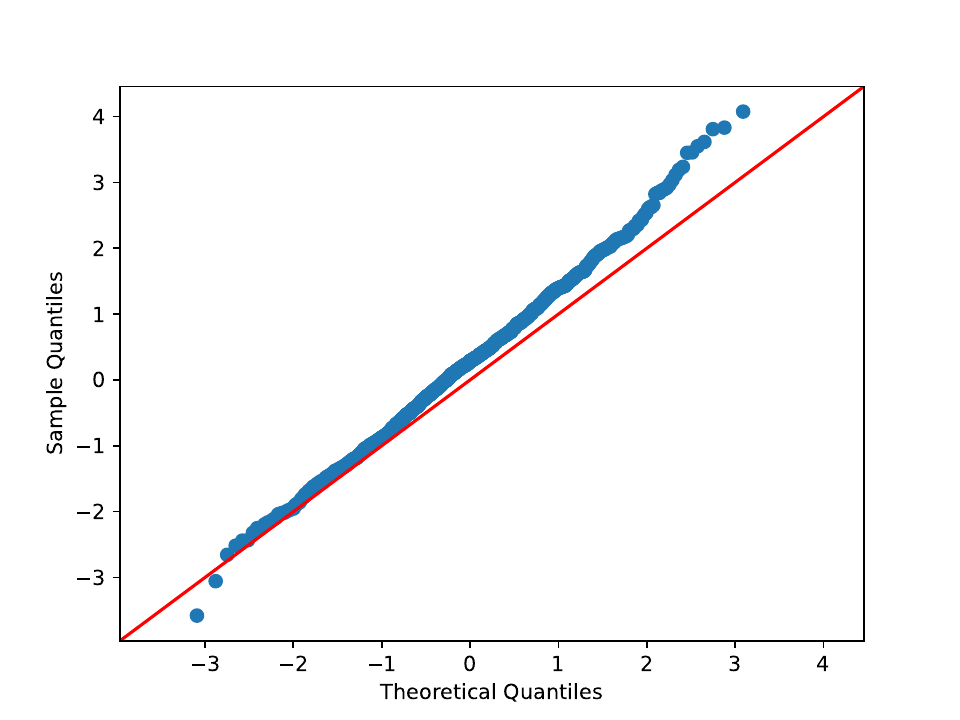}
           \caption{Normal QQ plot}
           \label{fig: normality multi 1}
       \end{subfigure}

    }
\caption{Asymptotic normality of episodic value functions}
\label{fig: normality multi}
\end{figure}

\begin{figure}[tb]
    {
        \centering
        \includegraphics[width=0.48\textwidth]{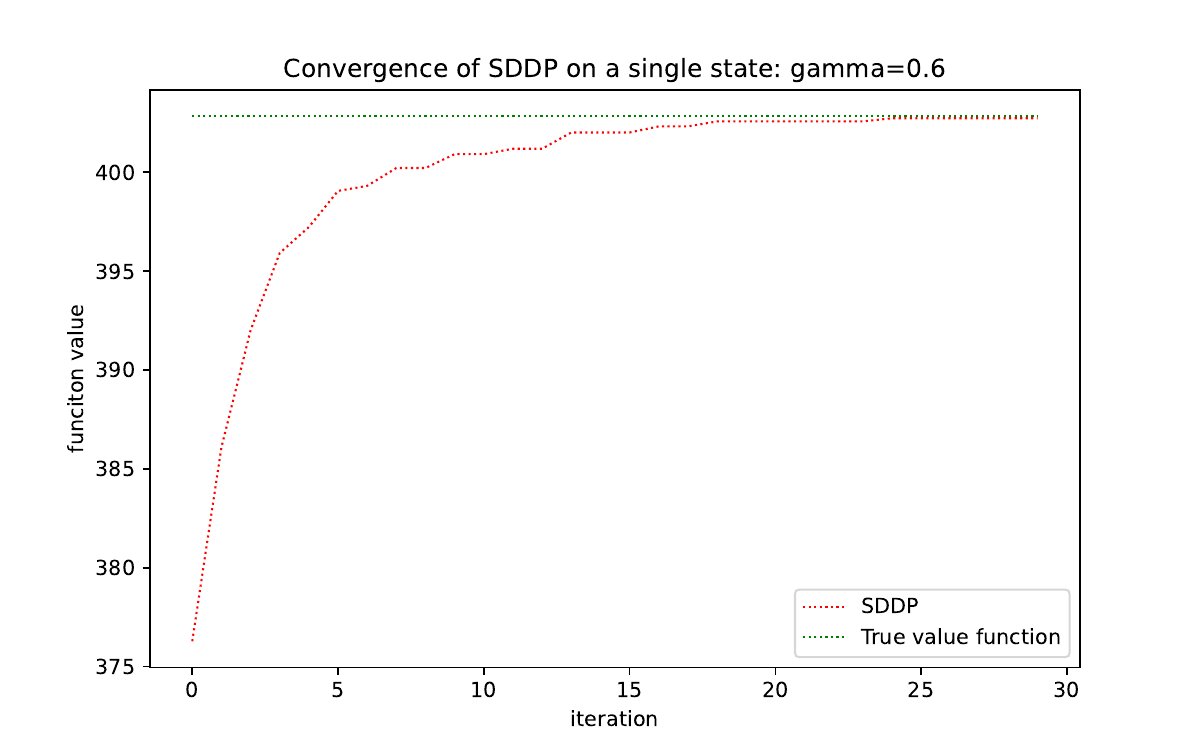}
        \includegraphics[width=0.48\textwidth]{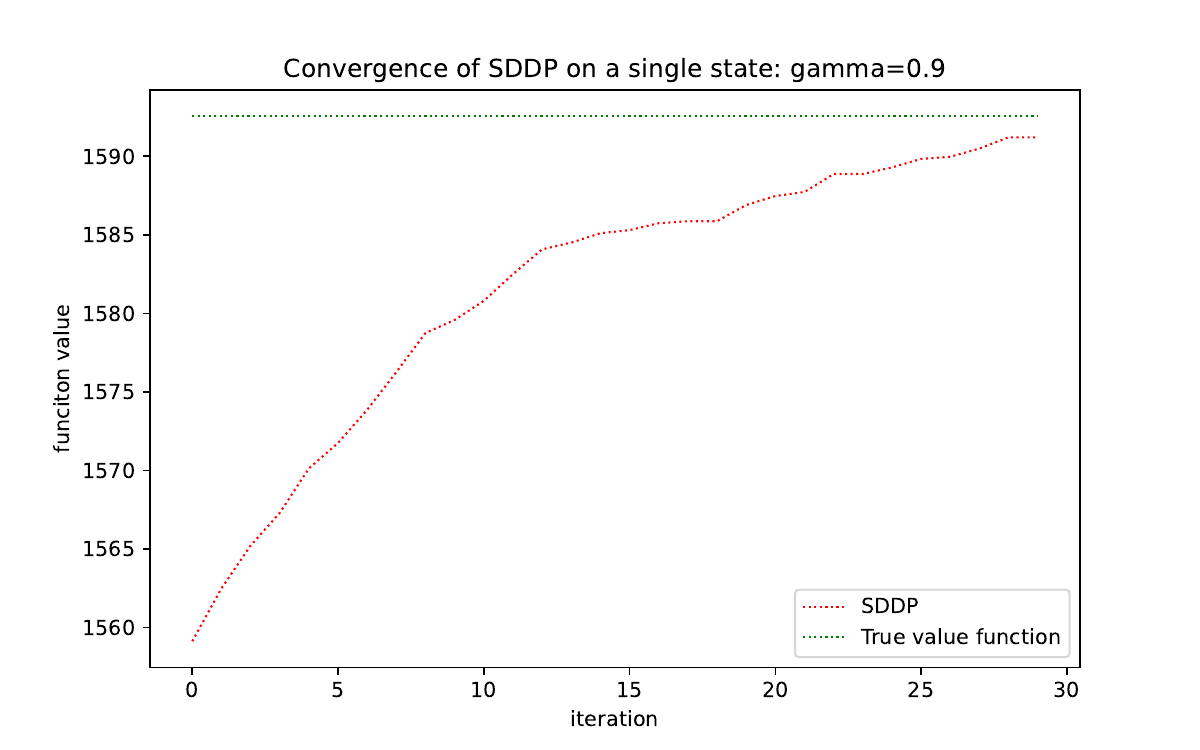}
    }
\caption{SDDP algorithm with single episode for 5-D Inventory Control.}
\label{fig: lower multi}
\end{figure}
\begin{figure}[tb]
    {
        \centering
        \includegraphics[width=0.48\textwidth]{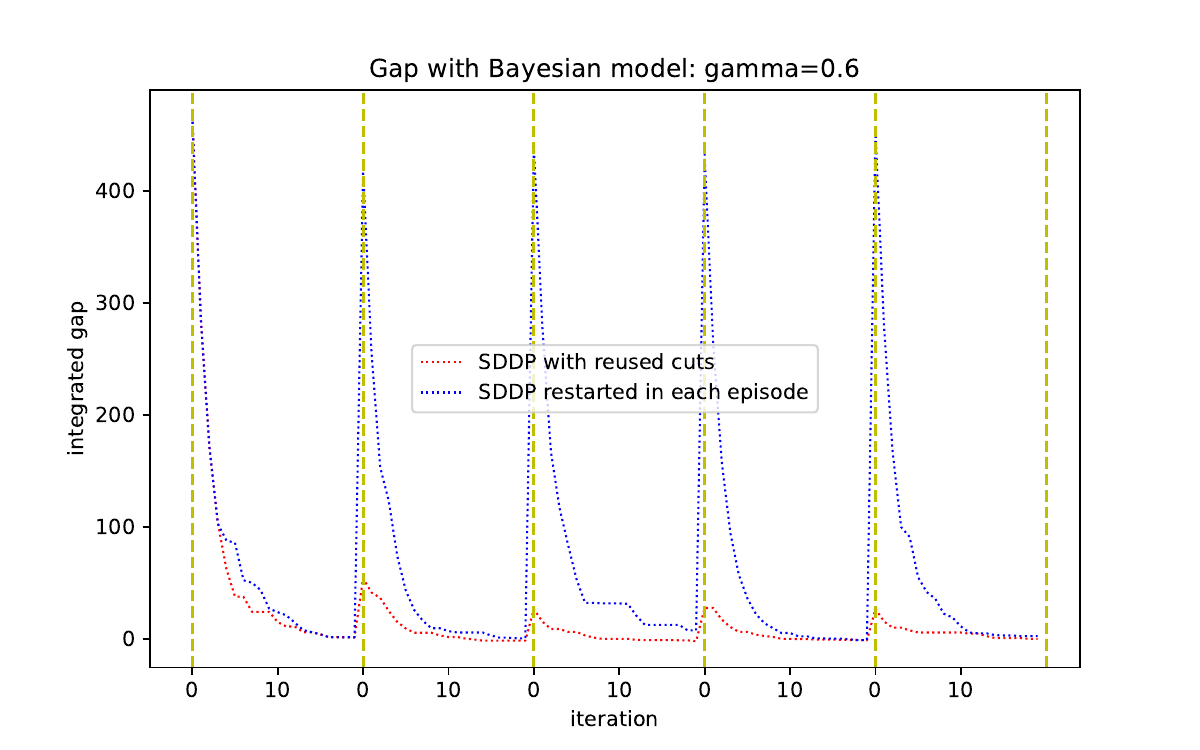}
        \includegraphics[width=0.48\textwidth]{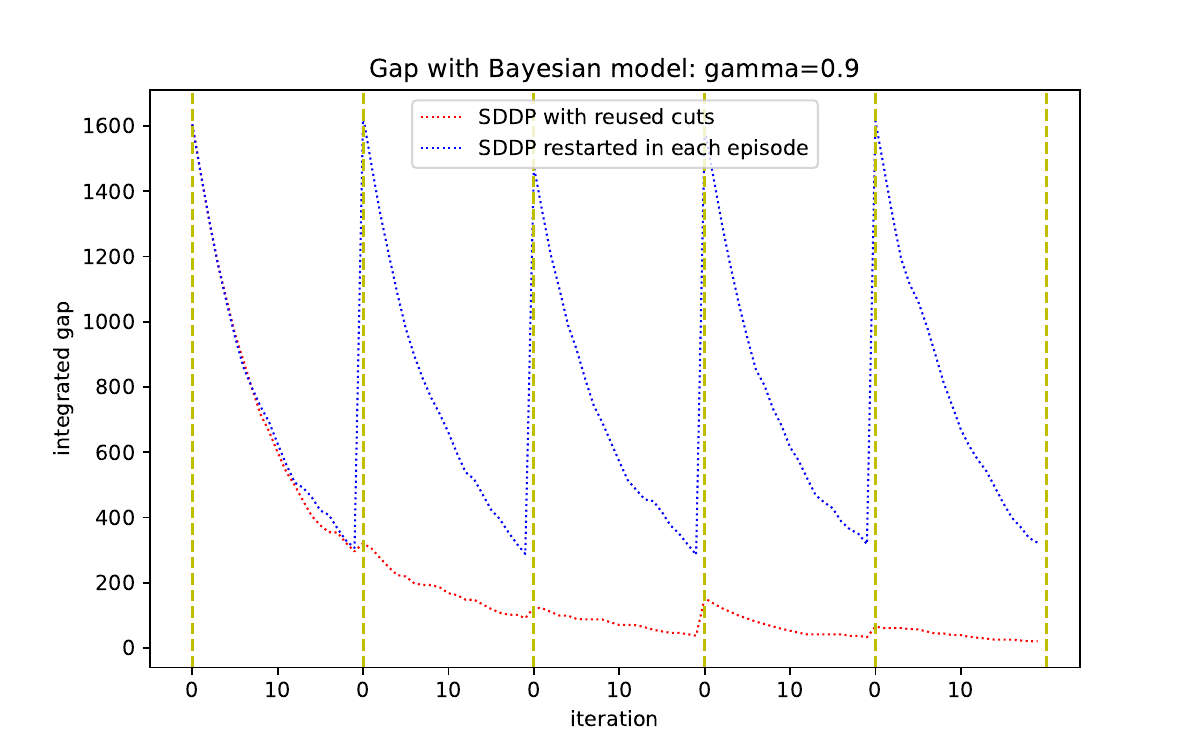}
    }
\caption{SDDP algorithm for the 5-D Inventory Control.}
\label{fig: 4}
\end{figure}

\subsection{Multi-dimensional case} \label{sec: numerical multi-dim}
We also test on the multi-dimensional version of this inventory example. Suppose we have $N\_prod = 5$ types of products.   The $i$th product has the order cost $c_i = 1+0.5*\sin{(i-1)}$, holding cost $h_i = 2+0.5*\sin{(i-1)} $ and backorder penalty cost $b_i = 3 + 0.5*\sin(i-1)$, for $i=1,2,\ldots,5$. The state is a five-dimensional vector, with each component being the current inventory level of each product. Assume the demand for each product follows an exponential distribution with unknown mean $\theta_i = 10 + 0.5 i, i=1,2,\ldots,5$.

We first show the convergence rate and asymptotic normality of the episodic value functions to the true optimal value function, and then demonstrate Algorithm~\ref{algorithm: episodic_SDDP} with or without reusing past lower approximation. Similar as for the 1-D case, due to the larger estimation error of $V^*_N$, in Figure \ref{fig: normality multi} we set the mean of demand distribution for the $i$th product to be $1+\sin(i-1)*0.2$ { and still use $10^5$ samples from $P_*$ or $P_N$ to approximate the expectation in $V^*$ and $V^*_N$. We evaluate $V^*_N$ with $50$ episodes for testing convergence in Figure \ref{fig: convergence multi}. In Figure \ref{fig: normality multi}, for testing asymptotic normality we evaluate $V^*_N$ at $N=100$ and at state $x_0 = \mathbf{0}$. The empirical pdf and QQ plot are calculated with 1000 replications. {In Figure \ref{fig: lower multi} and \ref{fig: 4}, we use sample size $M=100$ in Algorithm~\ref{algorithm: episodic_SDDP} and run for 5 episodes with number of SDDP iterations $K=30$ for each episode. In Figure \ref{fig: lower multi}, the lower approximation is evaluated at state $x_0 = \mathbf{0}$.  }

Similar results hold as in 1-D case in Figures \ref{fig: convergence multi}, \ref{fig: normality multi}, \ref{fig: lower multi} and \ref{fig: 4}. It is worth noting that there are more benefits of reusing past lower approximation in the higher-dimensional problem, in particular, because SDDP converges slower in high dimensions. For example, in the second plots of Figures~\ref{fig: lower multi} and \ref{fig: 4} with discount factor $\gamma = 0.9$, the difference between the two algorithms (with and without warm start) is starker, and we note that with warm start the percentage of reused previous cuts is higher than $80\%$ in the last three episodes.

\subsection{Performance Comparison}
In Section \ref{sec: numerical one-dim} and \ref{sec: numerical multi-dim}, we numerically show the asymptotic normality of $V^*_N$ as well as the convergence of the SDDP algorithm. In this section, we compare the proposed episodic Bayesian optimal control (EBOC)  approach with
% with two existing methods. 
% {\color{red} (Yuhao: Not sure about the validity of this sentence.) To our best knowledge, there is no existing literature in optimal control considering the discounted infinite-horizon criteria and the episodic setting as this paper. } To make the comparison,  we employ 
an extension of the distributionally robust stochastic control (denoted by DRSC), proposed in \cite{shapiro2012minimax}, and the LazyPSRL method, a variant of posterior sampling for online reinforcement learning (RL) from \cite{osband2013more,abbasi2015bayesian}.

The DRSC algorithm is an extension of the algorithm in \cite{shapiro2012minimax}, which considered a minimax formulation with the (robust) value function defined as 
$$
V_{\text {rob }}(x)=\inf _{u \in \mathcal{U}} \sup _{Q \in \mathcal{M}} \mathbb{E}_{\xi\sim Q}\left[c(x, u, \xi)+\gamma V_{\text {rob }}(F(x, u, \xi))\right],
$$
where $\xi\sim Q$ and $\mathcal{M}$ is he (ambiguity) set of probability distributions $Q$. Notably, \cite{shapiro2012minimax} only considers a single episode, meaning that the ambiguity set $\mathcal{M}$ is only constructed once with some pre-collected historical data set.  
To apply this to the episodic setting, we extend that approach by using all data collected from the previous episodes at the beginning of each episode and shrinking the ambiguity set in an order of $O(\frac{1}{\sqrt{t}})$ (where $t$ is the current iteration). The rate  $O(t^{-1/2})$ is motivated by the asymptotics of the distributionally robust optimization in the static setting (cf. \cite{blanchet2023statistical}). As for the LazyPSRL algorithm, it computes a policy by first sampling an RL environment (in our context, $\theta$) and then solve for a (near) optimal policy in that sampled environment. In contrast, we generate multiple samples of $\theta$ to approximate the Bayesian average objective and compute the policy by solving the approximated problem. 

To evaluate the performance of an algorithm, we define the regret at iteration $t$ as
$$\operatorname{regret}(t) = V^{\pi_t}(x_t) - V^*(x_t),$$
where $\pi_t$ is the policy computed from the algorithm, $x_t$ is the state visited at iteration $t$, $V^{\pi_t}(x_t)$ is the value function at state $x_t$ with policy $\pi_t$ under the true parameter $\theta^*$, and $V^*(x_t)$ is the optimal value function at state $x_t$ under the true parameter $\theta^*$.

For fair comparison,  we want to make the optimization error (i.e., the error between the computed policy and the optimal policy) and evaluation error (computation error of $V^*$ and $V^{\pi_t})$) for different methods extremely small to highlight the performance difference due to the different formulations by different methods. To achieve this, we test on the same single-product inventory example in Section \ref{sec:inventory}  with the unknown demand distribution to be a Poisson distribution. More details on the implementation of the algorithms (i.e., how to compute $\pi_t$ for different algorithms and how to evaluate the policy) are provided in the electronic companion \ref{EC: implementation}.

On a related note, as mentioned in Section \ref{sec: intro}, \cite{duff2002optimal} and \cite{Zhou2022NIPS} employ the Bayes adaptive MDP and Bayesian risk MDP, respectively, both of which treat the posterior distribution as an augmented state which varies over horizons. While their models take into account all possible updating of the posterior distribution, computationally solving the resulted MDPs  suffer from the curse of dimensionality due to the continuous belief state (i.e., the posterior distribution), even when the original state space and action space are finite. In particular, to solve the Bayes adaptive (Bayesian risk) MDP, both \cite{duff2002optimal} and \cite{Zhou2022NIPS} use dynamic programming, assuming a finite decision horizons. The computation cost increases exponentially in the number of total horizons, making their algorithms inapplicable to the scenarios with infinite horizon as studied in this paper. % To facilitate adaptions of their approach to infinite horizons, more computational efficient algorithm needs to be developed, which is beyond the scope of this paper.

\subsubsection*{Experiment details.}
In the inventory control problem, the unit order cost $c=1$, holding cost $h=2$, backorder penalty cost $b=3$, the demand follows a Poisson distribution with unknown mean $5$, the discount factor $\gamma = 0.9$, the prior distribution for all algorithms to be a gamma distribution $\Gamma(1,1)$, pre-collected initial data size $10$. For the proposed EBOC, we run with two different sample sizes $2,5$ for $\theta$. For the episodic length (number of iterations within a episode), we first test with a constant episodic length $H=5$ (see Figure \ref{fig: comparison constant}) and then test with the episodic length determined by the LazyPSRL (see Figure \ref{fig: comparison lazy}), which re-calculates a new policy only when the variance of the posterior distribution is reduced sufficiently. In both figures we plot average regret as well as its $95\%$ confidence interval (shown in shaded bands) by running 500 macro-replications.
\begin{figure}[h]
\includegraphics[width=\linewidth]{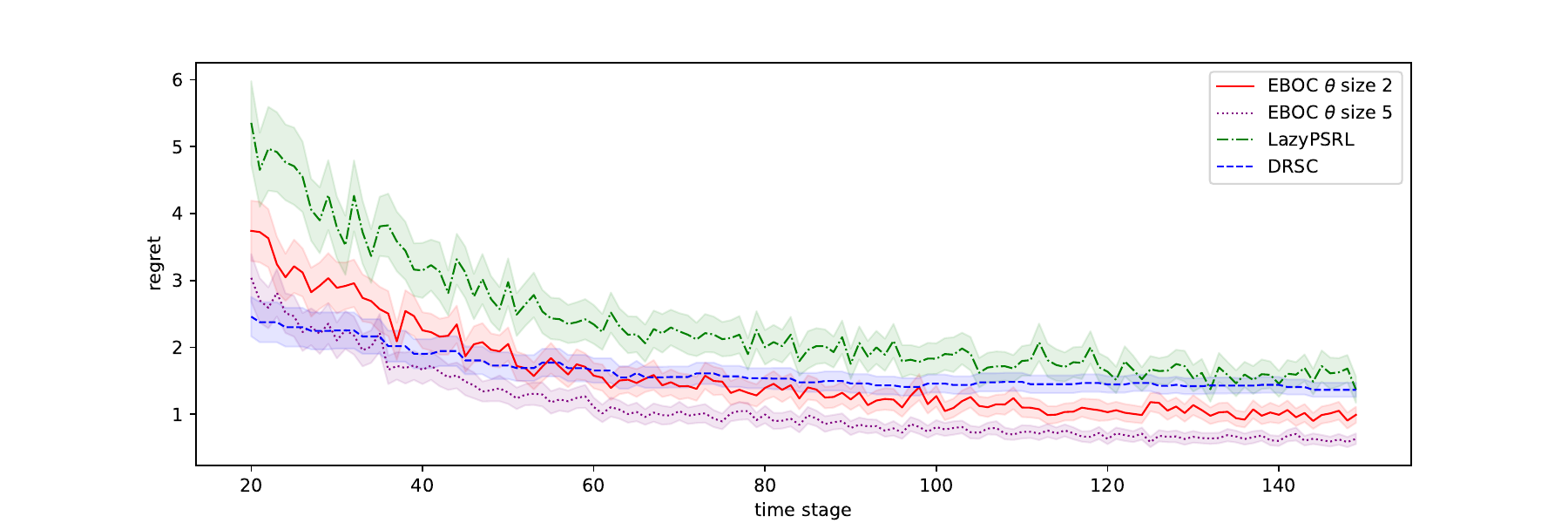}
    \caption{Regret over time: constant episodic length.}
    \label{fig: comparison constant}
\end{figure}
\begin{figure}[h]
\includegraphics[width=\linewidth]{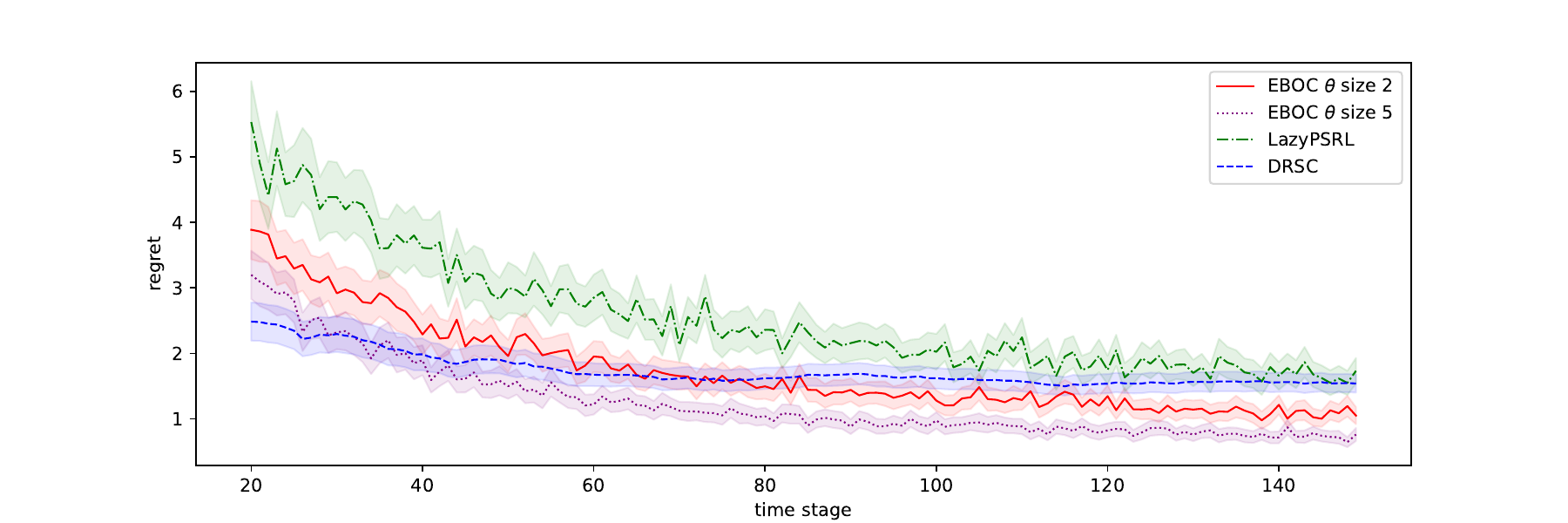}
    \caption{Regret over time: episodic length by lazyPSRL}
    \label{fig: comparison lazy}
\end{figure}

\subsubsection*{Result:}
\begin{enumerate}
    \item In both Figure \ref{fig: comparison constant} and \ref{fig: comparison lazy}, the proposed EBOC framework with $\theta$ sample size $5$ performs the best and EBOC with $\theta$ sample size $2$ performs the second best, showing the efficiency of EBOC.
    
    \item EBOC outperforms LazyPSRL for the following reason. In LazyPSRL, the posterior distribution is used for posterior sampling, which helps balance the trade-off between exploitation and exploration in RL. However, in the control setting, the random data $\xi_t$ is independent of the control policy, indicating that there is no benefit of exploration since data collected by any policy contains the same information. Consequently, the regret of LazyPSRL grows as the difference between  the sampled environment $\theta$ by the LazyPSRL and the true environment $\theta^*$. In contrast, EBOC uses the entire posterior distribution to aggregate the model estimation and  solves the Bayesian average objective, as opposed to a sampled objective in LazyPSRL. The Bayesian average counterpart of the Bellman equation \eqref{bellm-3} can be regarded as an estimation of the true optimal value function, which fully exploits the current information. The sample size of $\theta$ is then chosen to control the SAA error rather than balancing the exploitation-exploration trade-off, which does not exist in the control problem. Indeed, as indicated by both Figure \ref{fig: comparison constant} and \ref{fig: comparison lazy}, with increased sample size, both the average regret and the variation are reduced.
    
    \item Compared with DRSC, EBOC performs better in the sense that the average regret of EBOC decreases much faster. We conjecture that the reason is that  the ambiguity set of DRSC shrinks too slowly, resulting in an overly conservative policy. However, to the best of our knowledge, there is no existing study on how to adaptively shrink the ambiguity set for distributionally robust approaches in the dynamic setting, a topic beyond the scope of this paper. Furthermore, the construction of a good ambiguity set may depend on the specific problem and require substantially more effort than the natural way of updating posteriors in EBOC.
\end{enumerate}

}
\section{Conclusions and Future Research}
In this paper, we propose an episodic Bayesian approach to address unknown randomness distribution in stochastic optimal control problems. It incorporates Bayesian learning of the unknown distribution with optimization of the control policy. In each episode, the current policy is exercised and more data on the randomness are collected to update the Bayesian posterior, which is in turn used to update the control policy in the next episode. We show, with some huristics, that the episodic value functions converge almost surely to the true optimal value function at an asymptotic rate of $1/\sqrt{N}$, where $N$ is the number of data points (which is also the number of episodes when the episode length is $1$). For a class of problems with convex cost functions and linear state dynamics (so that the value functions are convex), we develop a computational algorithm based on the stochastic dual dynamic programming (SDDP) approach, with a warm start for each episode by reusing the value function approximation from the previous episode. Our numerical results on an inventory control problem verify the theoretical results and also demonstrate the effectiveness of the proposed algorithm. Our approach also numerically outperforms LazyPSRL, a variant of posterior sampling approach for reinforcement learning, and a distributionally robust approach with shrinking ambiguity set, in terms of the regret of the value function.

There are several directions for future research. First, we have assumed that the length of each episode is constant, but it could be time-varying or even random in accordance with the Bayesian learning of the unknown distribution. For example, the episode length could become longer as the Bayesian posterior stabilizes at later times. We have also assumed the Bayesian average problem is solved to optimality in each episode, but it is not always possible or necessary to do so. There should be a balance between the estimation error (due to Bayesian learning from finite data) and optimization error (due to sub-optimality of the solution to the episodic Bayesian average problem). Second, as we mentioned in the paper, the assumed parametric form for the unknown randomness distribution could lead to model misspecification, and therefore, the Bayesian distributionally robust optimization approach for static (one-stage) optimization problems (as proposed in \cite{SZY2023})  could be a valuable direction to investigate for the dynamic problem here.

\section*{Acknowledgement}
All authors are grateful for the support by Air Force Office of Scientific Research (AFOSR)
under Grant FA9550-22-1-0244. The second, third, and fourth authors are also grateful for the support by National Science Foundation (NSF) under Grant DMS2053489 and NSF AI Institute for Advances in Optimization (AI4OPT).

%Bibliography
\bibliographystyle{plainnat}  
\bibliography{references}

\begin{thebibliography}{37}
\providecommand{\natexlab}[1]{#1}
\providecommand{\url}[1]{\texttt{#1}}
\expandafter\ifx\csname urlstyle\endcsname\relax
  \providecommand{\doi}[1]{doi: #1}\else
  \providecommand{\doi}{doi: \begingroup \urlstyle{rm}\Url}\fi

\bibitem[Abbasi-Yadkori and Szepesv{\'a}ri(2015)]{abbasi2015bayesian}
Yasin Abbasi-Yadkori and Csaba Szepesv{\'a}ri.
\newblock Bayesian optimal control of smoothly parameterized systems.
\newblock In \emph{UAI}, pages 1--11, 2015.

\bibitem[Abeille and Lazaric(2018)]{abeille2018improved}
Marc Abeille and Alessandro Lazaric.
\newblock Improved regret bounds for thompson sampling in linear quadratic control problems.
\newblock In \emph{International Conference on Machine Learning}, pages 1--9. PMLR, 2018.

\bibitem[Basar and Bernhard(2008)]{Basar2008Hinfty}
T.~Basar and P.~Bernhard.
\newblock \emph{$H^\infty$-optimal control and related minimax design problems -- A dynamic game approach}.
\newblock Birkhauser Boston, 2008.

\bibitem[Bertsekas and Shreve(1978)]{ber78}
D.P. Bertsekas and S.E. Shreve.
\newblock \emph{Stochastic Optimal Control, The Discrete Time Case}.
\newblock Academic Press, New York, 1978.

\bibitem[Bielecki et~al.(2019)Bielecki, Chen, Cialenco, Cousin, and Jeanblanc]{Bielecki2019ARC}
Tomasz~R. Bielecki, Tao Chen, Igor Cialenco, Areski Cousin, and Monique Jeanblanc.
\newblock Adaptive robust control under model uncertainty.
\newblock \emph{SIAM Journal on Control and Optimization}, 57\penalty0 (2):\penalty0 925--946, 2019.

\bibitem[Billingsley(1995)]{Billingsley}
P.~Billingsley.
\newblock \emph{Probability and Measure}.
\newblock John Wiley and Sons, New York, 1995.

\bibitem[Blanchet and Shapiro(2023)]{blanchet2023statistical}
Jose Blanchet and Alexander Shapiro.
\newblock Statistical limit theorems in distributionally robust optimization.
\newblock In \emph{2023 Winter Simulation Conference (WSC)}, pages 31--45. IEEE, 2023.

\bibitem[Doob(1948)]{doob1948application}
Joseph~L. Doob.
\newblock Application of the theory of martingales.
\newblock \emph{Actes du Colloque International Le Calcul des Probabilites et ses applications}, pages 23--27, 1948.

\bibitem[Duff(2002)]{duff2002optimal}
Michael~O'Gordon Duff.
\newblock \emph{Optimal Learning: Computational procedures for Bayes-adaptive Markov decision processes}.
\newblock University of Massachusetts Amherst, 2002.

\bibitem[Fernholz(1983)]{Fernholz1983vonMises}
L.~T. Fernholz.
\newblock \emph{von Mises Calculus for Statistical Functional, Lecture Notes in Statistics, No. 19}.
\newblock New York: Springer-Verlag, 1983.

\bibitem[Gilboa and Schmeidler(1989)]{Gilboa1989MaxminEU}
Itzhak Gilboa and David Schmeidler.
\newblock Maxmin expected utility with non-unique prior.
\newblock \emph{Journal of Mathematical Economics}, 18:\penalty0 141--153, 1989.

\bibitem[Gonz\'{a}lez-Trejo et~al.(2002)Gonz\'{a}lez-Trejo, Hern\'{a}ndez-Lerma, and Hoyos-Reyes]{Gonz2002MinimaxControl}
J.~I. Gonz\'{a}lez-Trejo, O.~Hern\'{a}ndez-Lerma, and L.~F. Hoyos-Reyes.
\newblock Minimax control of discrete-time stochastic systems.
\newblock \emph{SIAM Journal on Control and Optimization}, 41\penalty0 (5):\penalty0 1626--1659, 2002.

\bibitem[Hansen et~al.(2006)Hansen, Sargent, Turmuhambetova, and Williams]{Hansen2006robust}
L.~P. Hansen, G.~Sargent, G.~Turmuhambetova, and N.~Williams.
\newblock Robust control and model misspecification.
\newblock \emph{J. Econom. Theory}, 2006.

\bibitem[Kumar and Varaiya(2015)]{Kumar2015BayesianAC}
P.~R. Kumar and P.~Varaiya.
\newblock \emph{Stochastic Systems: Estimation, Identification, and Adaptive Control}, chapter Chapter 11: Bayesian Adaptive Control.
\newblock Prentice Hall, 2015.

\bibitem[Lan and Shapiro(2024)]{LanShapiro:2024survey}
Guanghui Lan and Alexander Shapiro.
\newblock Numerical methods for convex multistage stochastic optimization.
\newblock \emph{Foundations and Trends in Optimization}, 6\penalty0 (2):\penalty0 63 -- 144, 2024.

\bibitem[Lim et~al.(2006)Lim, Shanthikumar, and Shen]{Lim2006Robust}
A.~E.~B. Lim, G.~J. Shanthikumar, and Z.~J.~Max Shen.
\newblock Model uncertainty, robust optimization and learning.
\newblock \emph{Tutorials in Operations Research, INFORMS}, pages 66--94, 2006.

\bibitem[Lin et~al.(2022)Lin, Ren, and Zhou]{Zhou2022NIPS}
Yifan Lin, Yuxuan Ren, and Enlu Zhou.
\newblock Bayesian risk {M}arkov decision processes.
\newblock In S.~Koyejo, S.~Mohamed, A.~Agarwal, D.~Belgrave, K.~Cho, and A.~Oh, editors, \emph{Advances in Neural Information Processing Systems}, volume~35, pages 17430--17442, 2022.

\bibitem[Liu et~al.(2024)Liu, Lin, and Zhou]{ZhouLiu2022bayesianSGD}
Tianyi Liu, Yifan Lin, and Enlu Zhou.
\newblock Bayesian stochastic gradient descent for stochastic optimization with streaming input data.
\newblock \emph{SIAM Journal on Optimization}, \penalty0 (1):\penalty0 389--418, 2024.

\bibitem[Osband and Van~Roy(2014)]{osband2014near}
Ian Osband and Benjamin Van~Roy.
\newblock Near-optimal reinforcement learning in factored mdps.
\newblock \emph{Advances in Neural Information Processing Systems}, 27, 2014.

\bibitem[Osband and Van~Roy(2017)]{osband2017posterior}
Ian Osband and Benjamin Van~Roy.
\newblock Why is posterior sampling better than optimism for reinforcement learning?
\newblock In \emph{International conference on machine learning}, pages 2701--2710. PMLR, 2017.

\bibitem[Osband et~al.(2013)Osband, Russo, and Van~Roy]{osband2013more}
Ian Osband, Daniel Russo, and Benjamin Van~Roy.
\newblock (more) efficient reinforcement learning via posterior sampling.
\newblock \emph{Advances in Neural Information Processing Systems}, 26, 2013.

\bibitem[Pereira and Pinto(1991)]{per1991}
M.V.F. Pereira and L.M.V.G. Pinto.
\newblock Multi-stage stochastic optimization applied to energy planning.
\newblock \emph{Mathematical programming}, 52\penalty0 (1-3):\penalty0 359--375, 1991.

\bibitem[Rieder(1975)]{Rieder1975BayesianDP}
U.~Rieder.
\newblock Bayesian dynamic programming.
\newblock \emph{Adv. Appl. Prob.}, pages 330--348, 1975.

\bibitem[Schwartz(1965)]{Schwartz1965OnBP}
Lorraine Schwartz.
\newblock On {B}ayes procedures.
\newblock \emph{Zeitschrift f{\"u}r Wahrscheinlichkeitstheorie und Verwandte Gebiete}, 4:\penalty0 10--26, 1965.

\bibitem[Shapiro and Cheng(2021)]{shacheng2021}
A.~Shapiro and Y.~Cheng.
\newblock Central limit theorem and sample complexity of stationary stochastic programs.
\newblock \emph{Operations Research Letters}, 49:\penalty0 676--681, 2021.

\bibitem[Shapiro et~al.(2021)Shapiro, Dentcheva, and Ruszczy\'{n}ski]{SDR}
A.~Shapiro, D.~Dentcheva, and A.~Ruszczy\'{n}ski.
\newblock \emph{Lectures on Stochastic Programming: Modeling and Theory}.
\newblock SIAM, Philadelphia, third edition, 2021.

\bibitem[Shapiro(2012)]{shapiro2012minimax}
Alexander Shapiro.
\newblock Minimax and risk averse multistage stochastic programming.
\newblock \emph{European Journal of Operational Research}, 219\penalty0 (3):\penalty0 719--726, 2012.

\bibitem[Shapiro et~al.(2023)Shapiro, Zhou, and Lin]{SZY2023}
Alexander Shapiro, Enlu Zhou, and Yifan Lin.
\newblock Bayesian distributionally robust optimization.
\newblock \emph{SIAM Journal on Optimization}, 33\penalty0 (2):\penalty0 1279--1304, 2023.

\bibitem[S{\^i}rbu(2014)]{Srbu2014ANO}
Mihai S{\^i}rbu.
\newblock A note on the strong formulation of stochastic control problems with model uncertainty.
\newblock \emph{Electronic Communications in Probability}, 19:\penalty0 1--10, 2014.

\bibitem[Strens(2000)]{strens2000bayesian}
Malcolm Strens.
\newblock A bayesian framework for reinforcement learning.
\newblock In \emph{ICML}, volume 2000, pages 943--950, 2000.

\bibitem[Theocharous et~al.(2017)Theocharous, Wen, Abbasi-Yadkori, and Vlassis]{theocharous2017posterior}
Georgios Theocharous, Zheng Wen, Yasin Abbasi-Yadkori, and Nikos Vlassis.
\newblock Posterior sampling for large scale reinforcement learning.
\newblock \emph{arXiv preprint arXiv:1711.07979}, 2017.

\bibitem[Tzortzis et~al.(2019)Tzortzis, Charalambous, and Charalambous]{Tzortzis2019TVA}
I.~Tzortzis, C.~D. Charalambous, and T.~Charalambous.
\newblock Infinite horizon average cost dynamic programming subject to total variation distance ambiguity.
\newblock \emph{SIAM Journal of Control and Optimization}, 57\penalty0 (4):\penalty0 2843--2872, 2019.

\bibitem[van~der Vaart(1998)]{vaart}
A.W. van~der Vaart.
\newblock \emph{Asymptotic Statistics}.
\newblock Cambridge Univ. Press, 1998.

\bibitem[Van~Parys et~al.(2016)Van~Parys, Kuhn, Goulart, and Morari]{Parys2016DRC}
Bart P.~G. Van~Parys, Daniel Kuhn, Paul~J. Goulart, and Manfred Morari.
\newblock Distributionally robust control of constrained stochastic systems.
\newblock \emph{IEEE Transactions on Automatic Control}, 61\penalty0 (2):\penalty0 430--442, 2016.

\bibitem[{Wu} et~al.(2018){Wu}, {Zhu}, and {Zhou}]{Wu:2018BRO}
D.~{Wu}, H.~{Zhu}, and E.~{Zhou}.
\newblock A {Bayesian} risk approach to data-driven stochastic optimization: Formulations and asymptotics.
\newblock \emph{SIAM Journal on Optimization}, 28\penalty0 (2):\penalty0 1588--1612, 2018.

\bibitem[Yang(2018)]{Yang2018WassersteinDR}
Insoon Yang.
\newblock Wasserstein distributionally robust stochastic control: A data-driven approach.
\newblock \emph{IEEE Transactions on Automatic Control}, 66:\penalty0 3863--3870, 2018.

\bibitem[Zipkin(2000)]{zipkin}
P.H. Zipkin.
\newblock \emph{Foundation of inventory management}.
\newblock McGraw-Hill, 2000.

\end{thebibliography}

\newpage
\appendix
\section{Discussion on state-dependent control} \label{EC: state-dependent control}
Suppose now that the control set depends on the state vector and is given in the form
\begin{equation}\label{appen-1}
 \U(x):=\{u\in \U:g_k(x,u)\le 0,\;k=1,...,K\}.
\end{equation}
The counterpart of the Bellman equation \eqref{bellm-1} is
 \begin{equation}
 \label{appen-2}
V(x)=\inf\left\{ \bbe
   \left [c(x,u,\xi)+\gamma
V \big(F(x,u,\xi) \big)\right]: u\in \U,\;g(x,u)\le 0\right\},
\end{equation}
where $g(x,u)=(g_1(x,u),...,g_K(x,u))$.
The counterparts of the value functions $V^*$ and $V^*_N$ are given as solutions of the  Bellman equation of the form \eqref{appen-2}
when the expectation is taken with respect to the corresponding  distributions  $f(\cdot|\theta^*)$ and $P_N$.
The value functions are convex if in addition to Assumption \ref{ass-conv} the functions $g_k(x,u)$, $i=1,...,K$, are convex in $(x,u)$.

For a current approximation $\underline{\V}_N$ of the value function the required subgradient at a trial point $\bar{x}$ is given by an optimal solution of the dual of the problem
\begin{equation}
 \label{appen-3}
\begin{array}{cll}
 \min\limits_{u\in \U,\,z\in \bbr^M}&\sum_{j=1}^M   c_j(\bar{x},u)+z_j \\
 {\rm s.t.} & \underline{\V}_N(A_j\bar{x}+B_j u+b_j)\le z_j,\;j=1,...,M,\\
 & g_k(\bar{x},u)\le 0,\;k=1,...,K.
\end{array}
 \end{equation}
If the set $\U$ is polyhedral and the  functions $\underline{\V}_N(\cdot)$, $c_j(\cdot,\cdot)$, $j=1,...,M$,  and $g_k(\cdot,\cdot)$, $k=1,...,K$,  are given as   maxima of the respective (finite) sets  of affine functions,   problem \eqref{appen-3} can be formulated as a linear programming problem.

\section{Uniform convergence of $v_N(x)$} \label{sec:EClemma}
\begin{lemma}
 \label{lem-uniform}
Consider bounded  functions $\phi(x,u), \phi_N(x,u):\X\times \U\to \bbr$. Suppose  that
$\phi_N(x,u)$ converges to $\phi(x,u)$ uniformly in $(x,u)\in \X\times \U$.  Then $\nu_N(x):=\inf_{u\in \U} \phi_N(x,u)$  converges to $\nu(x):=\inf_{u\in \U} \phi(x,u)$   uniformly in $x\in \X$.
\end{lemma}

{\bf Proof.} Let
\[
\kappa_N:=\sup_{(x,u)\in \X\times \U}|\phi_N(x,u)-\phi(x,u)|.
\]
By the assumption, $\kappa_N$ tends to zero as $N\to \infty$.
We have that
\[
\sup_{x\in \X}\{\nu_N(x)-\nu(x)\}\le \sup_{(x,u)\in \X\times \U}\{\phi_N(x,u)-\phi(x,u)\}\le \kappa_N,
\]
and
\[
\sup_{x\in \X}\{\nu(x)-\nu_N(x)\}\le \sup_{(x,u)\in \X\times \U}\{\phi(x,u)-\phi_N(x,u)\}\le \kappa_N,
\]
and hence the conclusion follows.  $\hfill\square$

\section{Danskin theorem}
\label{danskin}
Formula \eqref{direc} for the directional derivative of the value function  is motivated by the so-called  Danskin theorem.
That is, let $\Omega$ be a  (Hausdorff) topological  space, $X$ be a Banach space,  $\phi(x,\omega)$ be a real valued function and
$\psi(x):=\inf_{\omega\in \Omega} \phi(x,\omega)$ be the respective optimal value function. Suppose that: $\phi(\cdot,\omega)$ is differentiable, $\phi(x,\omega)$ and  the derivative $D_x\phi(x,\omega)$ are continuous on $X\times \Omega$,  and  $\Omega$ is compact. Then $\psi(\cdot)$ is directionally differentiable with the directional derivative
\[
\psi'(x,d)=\inf_{\omega\in \Omega^*(x)} D_x \phi(x,\omega)d,
\]
where $ \Omega^*(x):=\arg\min_{\omega\in \Omega} \phi(x,\omega)$  (e.g., Theorem 4.13 in Bonnans and Shapiro (2000)). In particular if  $\Omega^*(x)=\{\bar{\omega}(x)\}$ is the singleton, then $\psi(\cdot)$ is
(Fr\'echet) differentiable
and $\nabla\psi(x)=\nabla_x\phi(x,\omega^*)$, where $\omega^*=\bar{\omega}(x)$.
The compactness assumption is needed to ensue  continuity of the minimizer function $\bar{\omega}(\cdot)$.
The technical difficulty in applying this to derivation of the directional derivative in \eqref{direc}  is that the optimal policy lives in an infinite dimensional space and the   assumption of compactness could be difficult to verify. This is discussed in  Shapiro and Cheng (2021). %Because of that this result is partially heuristic. Not much we can do about that. It appears that this works OK in applications.

\section{Interchangeability of limit and expectation}
\label{sec:ECinterchange}
The standard condition for such  interchangeability   is the following so-called uniform integrability:
\begin{equation}\label{unifint}
 \lim_{c\to\infty} \sup_{N\in \bbn}\int_{\{|W_N|\ge c\}}|W_N(\xi,\theta)|d \theta  d\xi= 0.
\end{equation}
In particular, \eqref{unifint} holds if a.s. $|W_N(\xi,\theta)|$ is dominated by an integrable function
(e.g.,  \cite[Theorems 9.36 and 9.37]{SDR}). This follows from the Lebesgue Dominated Convergence Theorem \cite{Billingsley}. A special case of this is that $|W_N(\xi,\theta)| \leq W$, where $W$ is a finite constant.

\section{Convexity of the value function  when the control set depends on the state vector.} \label{EC:convexity}
Suppose   that the control set depends on the state vector and is given in the form
\begin{equation}\label{appen-1}
 \U(x):=\{u\in \U:g_k(x,u)\le 0,\;k=1,...,K\}.
\end{equation}
The counterpart of the Bellman equation \eqref{bellm-1} is
 \begin{equation}
 \label{appen-2}
V(x)=\inf\left\{ \bbe
   \left [c(x,u,\xi)+\gamma
V \big(F(x,u,\xi) \big)\right]: u\in \U,\;g(x,u)\le 0\right\},
\end{equation}
where $g(x,u)=(g_1(x,u),...,g_K(x,u))$.
The counterparts of the value functions $V^*$ and $V^*_N$ are given as solutions of the  Bellman equation of the form \eqref{appen-2}
when the expectation is taken with respect to the corresponding  distributions  $f(\cdot|\theta^*)$ and $P_N$.
The value functions are convex if in addition to Assumption \ref{ass-conv} the functions $g_k(x,u)$, $i=1,...,K$, are convex in $(x,u)$.

\section{Calculation of the subgradient.} \label{EC:subgradient}
For a current approximation $\underline{\V}_N$ of the value function the required subgradient at a trial point $\bar{x}$ is given by an optimal solution of the dual of the problem
\begin{equation}
 \label{appen-3}
\begin{array}{cll}
 \min\limits_{u\in \U,\,z\in \bbr^M}&\sum_{j=1}^M   c_j(\bar{x},u)+z_j \\
 {\rm s.t.} & \underline{\V}_N(A_j\bar{x}+B_j u+b_j)\le z_j,\;j=1,...,M,\\
 & g_k(\bar{x},u)\le 0,\;k=1,...,K.
\end{array}
 \end{equation}
If the set $\U$ is polyhedral and the  functions $\underline{\V}_N(\cdot)$, $c_j(\cdot,\cdot)$, $j=1,...,M$,  and $g_k(\cdot,\cdot)$, $k=1,...,K$,  are given as   maxima of the respective (finite) sets  of affine functions,   problem \eqref{appen-3} can be formulated as a linear programming problem.

\section{Numerical implementation detail} \label{EC: implementation}
For lazy PSRL, the optimal policy can be computed explicitly as in Section \ref{sec:inventory}. For EBOC, having sample $\theta_1,\theta_2,\ldots,\theta_n$, the optimal policy can be computed as the $\kappa = \frac{b-(1-\gamma)c}{b+h}$-quantile of a random variable $Y = \frac{1}{n}\sum_{i=1}^n Y_i$, where $Y_1,\ldots,Y_n$ are independent and $Y_i \sim \operatorname{Poisson}(\theta_i)$. We compute this quantile using Monte Carlo simulation with sample size $100$. For DRSC, we consider the ambiguity set with form
$$\mathcal{M} = \left\{Q: D_{\operatorname{KL}}(Q|| \operatorname{Poisson}(\widehat{\theta}_t)) \le \frac{1}{\sqrt{t}}\right\},$$
where $\widehat{\theta}_t$ is the sample average of all collected data until iteration $t$. Then, the distributionally robust optimal policy can be computed as $\pi_{rob}(x) = (x^*_{rob} -x )^+$ with $x^*_{rob} = \left\{
\begin{aligned}
   & \widehat{\theta}_t - \frac{1}{2\sqrt{t}} ~~ &\operatorname{ if } \kappa < \frac{1}{2}\\
   & \widehat{\theta}_t + \frac{1}{2\sqrt{t}} ~~& o.w.
\end{aligned}
\right.$.

To evaluate a policy, note all aforementioned optimal policies are order-up-to policies with some order-up-to level $x^*$. The value function can then be computed as (see \cite{shacheng2021}) 
$$V(x) = \left\{
\begin{aligned}
&
-c x+(1-\gamma)^{-1} \mathbb{E}_{\theta^*}\left[\gamma c D+(1-\gamma) c x^*+\psi\left(x^*, D\right)\right],
 ~~~~&\text {for } x \leq x^*,\\
&\mathbb{E}_{\theta^*}[\psi(x, D)+\gamma V(x-D)],~~~~& \text {for } x \geq x^* .
\end{aligned}
\right.$$
Since $D$ follows Poisson distribution, when $x > x^*$, we can compute
$$V(x) = (1-\gamma \mathrm{e}^{-\theta^*})^{-1} \sum_{k=1}^\infty \frac{(\theta^*)^k}{k!} \mathrm{e}^{-\theta^*} V (x-k).$$
And we truncate the infinite summation at $k=30$.

\end{document}